\documentclass[final, 12pt]{article}
\usepackage[english]{babel}
\usepackage[latin1]{inputenc}
\usepackage{hyperref}
\usepackage{amssymb,amsmath,amsthm}
\usepackage{enumerate}
\usepackage[conditional,light,first,bottomafter]{draftcopy}
\draftcopyName{DRAFT\space\today}{130}
\draftcopySetScale{65}
\usepackage[letterpaper, hmargin=2cm,vmargin={2.5cm,2.5cm}]{geometry}
\usepackage{setspace}
\makeatletter
\renewcommand{\section}{\@startsection%
{section}%
{1}%
{0em}%
{1.7em}%
{1.2em}%
{\normalfont\large\centering\bfseries}}
\renewcommand{\@seccntformat}[1]%
{\csname the#1\endcsname.\hspace{0.5em}}
\makeatother

\numberwithin{equation}{section}
\usepackage{fancyhdr}
\usepackage[usenames]{color}

\numberwithin{equation}{section}
\newtheorem{theorem}{Theorem}[section]

\newtheorem{proposition}{Proposition}[section]
\newtheorem{lemma}{Lemma}[section]
\newtheorem{corollary}{Corollary}[section]
\theoremstyle{definition}
\newtheorem{definition}{Definition}[section]

\theoremstyle{remark}
\newtheorem{remark}{Remark}[section]

\def\ocirc#1{\ifmmode\setbox0=\hbox{$#1$}\dimen0=\ht0 \advance\dimen0
  by1pt\rlap{\hbox to\wd0{\hss\raise\dimen0
  \hbox{\hskip.2em$\scriptscriptstyle\circ$}\hss}}#1\else {\accent"17 #1}\fi}
\DeclareMathOperator{\re}{Re}
\DeclareMathOperator{\im}{Im}
\DeclareMathOperator{\dom}{dom}

\DeclareMathOperator*{\res}{Res}

\begin{document}
\begin{titlepage}
\title{Spectral analysis of Volterra integrodifferential equations with the kernels, depending on parameter} 


\author{
\textbf{Romeo Perez Ortiz\footnote{Email: cemees.romeo@gmail.com}, Victor V. Vlasov\footnote{Email: vikmont@yandex.ru} and Nadezhda A. Rautian\footnote{Email: nrautian@mail.ru}}
\\[6mm]
\small Faculty of Mechanics and Mathematics \\[-1.6mm]
\small Lomonosov Moscow State University \\[-1.6mm]
\small Leninskie Gory, Moscow, 119991, Russian Federation\\[1mm]
}
\date{}
\maketitle
\vspace{4mm}
\begin{center}
\begin{minipage}{5in}
  \centerline{{\bf Abstract}} \bigskip Spectral analysis of operator-functions which are the symbols of the abstract integrodifferential equations of the Gurtin-Pipkin is provided. These equations represent abstract wave equations disturbed by terms involving Volterra operators. Correct solvability  in the Sobolev space $W_{2}^{2}((0, T), A^2)$, for arbitrary $T>0$, of that abstract integrodifferential equations is also studied. 
 \end{minipage}
\end{center}
\thispagestyle{empty}
\end{titlepage}
\section{Introduction}\label{introduction}
The paper is concerned with integrodifferential equations with unbounded operator coefficients in a~Hilbert space. 
The main part ($d^2u/dt^2+A^2u$) of the equation under consideration is an abstract hyperbolic-type equation disturbed by terms involving Volterra operators. These equations can be looked upon
as an abstract form of the Gurtin--\allowbreak Pipkin equation describing thermal phenomena and heat transfer in materials with memory or wave propagation in viscoelastic media.
A~complete analysis and  abundant examples of such equations in Banach and Hilbert spaces can be found in \cite{{GMJ}, {Dafermos}, {MFBL}, {MKR}, {MKR1}, {MRWR}, {MRWD}, {ESP}}.

Consider the following class of second-order abstract models
 \begin{align}
\frac{d^2u}{dt^2}+A^2u- &\int_{0}^{t}K(t-s)A^{2\theta} u(s) ds=f(t), \hspace{0.3cm} t \in \mathbb{R}_{+}, \label{system 1.1}\\
&u(+0)=\varphi_0, \hspace{0.4cm} u^{(1)}(+0)=\varphi_1 \label{system 1.2},
\end{align} where $A$ is a positive self-adjoint  operator with domain $\dom (A)\subset H$, $H$ is a Hilbert space. The variable  $\theta$ is a real number in $[0,1]$ and $K(t)$ is the kernel associated with the equation \eqref{system 1.1} (the Gurtin--\allowbreak Pipkin equation if $\theta=1$).  This type of equations appear  in various branches  of
mechanics and  physics, for instance, in heat transfer with finite propagation speed \cite{PG}, theory of viscoelastic media \cite{Dafermos},  kinetic theory of gases \cite{GK}, and thermal systems with memory \cite{Vegni}.
 
In this paper we restrict ourselves to two results: correct solvability of the initial value problem  \eqref{system 1.1}--\eqref{system 1.2} in the Sobolev space $W_{2}^{2}((0, T), A^2)$ for anyone $T>0$ and complete spectral analysis of operator-valued function $L(\lambda)$:  
\begin{align*}
L(\lambda)= \lambda^2I+A^2-\widehat{K}(\lambda)A^{2\theta},\hspace{0.2cm} \theta \in [0, 1],
\end{align*}which is the symbol of  integrodifferential equation \eqref{system 1.1}. The main goal of this paper is to study the location of spectrum of the operator-valued function $L(\lambda)$ in the case when $\theta \in [0, 1]$. Spectral analysis of integrodifferential equation \eqref{system 1.1}, in the case $\theta=1$,  was carried out in detail in the works \cite{VR, VR1,VR2, vlasovmrau2016, VRShamaev}.  The presence of a parameter $\theta \in [0, 1)$  changes significantly  the structure of the spectrum of the operator-valued function $L(\lambda)$ (see sections \ref{subsection-proof-of-theorem-about-spectral-structure-1}, \ref{subsection-proof-of-theorem-about-spectral-structure-when-condition-not-hold} and figures, given at the end of section \ref{subsection-proof-of-theorem-about-spectral-structure-when-condition-not-hold}). It is pertinent to mention that the equation \eqref{system 1.1} has been studied by many authors (see, for example, the monograph \cite{GMJ} and the bibliography given in it; see also works \cite{VR, VRShamaev, VR1,VR2, vlasovmrau2016, JMF, JME, JM} and the bibliography given in it). Results related to asymptotic behavior of solutions for systems with memory, for different $\theta \in [0, 1]$, were extensively studied in recent years (see  \cite{{MFBL}, {JMF}, {JME},{JM}} and the references given therein). In~\cite{JMF}, for instance, Mu\~noz Rivera and co-authors showed that the solutions for system \eqref{system 1.1}--\eqref{system 1.2} with $\theta \in [0, 1)$  decay polynomially as $t\to +\infty$, even if the kernel $K(t)$ decays exponentially.  In \cite{MFBL}, assuming the exponential decay of kernel $K(t)$ and $\theta=1$ Fabrizio and Lazzari proved the exponential decay of the solutions for system \eqref{system 1.1}--\eqref{system 1.2}.
In \cite{JME}, for the case $\theta=0$, Mu\~noz Rivera and co-authors showed that for the ionized atmosphere the dissipation produced by the conductivity kernel alone is not enough to
produce an exponential decay of the solution of an integrodifferential equation. In \cite{JM}, for the case $\theta=1$, Mu\~noz Rivera and Maria Naso proved that the solution of model \eqref{system 1.1}--\eqref{system 1.2} decays exponentially to zero if so does the kernel~$K(t)$.

In the section  \ref{subsection-about-distribution-of zeroes-meromorph-function-ell-n} is given a theorem including the distribution of spectrum of the operator-valued function $L(\lambda)$ on the left half-plane. From here the question  about the stability of solutions of abstract integro-differential equations of the form \eqref{system 1.1} when $\theta \in [0, 1]$ naturally arises. In this direction, there are several results related to asymptotic behavior of solutions for systems with memory, for different $\theta \in [0, 1]$ (see  \cite{{MFBL}, {JMF}, {JME},{JM}} and the references given therein). We note, however, that in the papers known to us, when $\theta \in [0, 1]$, and, in particular, in the works \cite{JMF, JME, JM} spectral analysis of the symbol $L(\lambda)$ of equation  \eqref{system 1.1} was not carry out. The results of this paper, therefore, are a natural development of the results of the works \cite{VR, VRShamaev, VR1,VR2, vlasovmrau2016}, in which spectral analysis was carried out in the case $\theta=1$. The information about the spectra of operator-valued function $L(\lambda)$ plays an important role in the analysis of behaviour of solutions of above mentioned equations. On this way it is possible to say wether or not the solutions will be decay exponentially. In the sections \ref{analysis-about-spectrum-in-the-easy-case} and \ref{representation-of-spectra-in-graphic-in-the-case-divergente-series} we formulate such results with explanations. Moreover, in the paper \cite{VR2}, when  $\theta=1$, was obtained the representation of solutions of the initial value problem \eqref{system 1.1}--\eqref{system 1.2} as series of exponential functions corresponding to spectrum of operator-valued function $L(\lambda)$. The formulation of similar results for an arbitrary $\theta \in [0, 1]$ was presented in the recent work \cite{RomeoRautian}.

 In the monograph \cite{GMJ} was proposed a semigroup approach in the researching initial value problem   \eqref{system 1.1}--\eqref{system 1.2}. To our mind the key moment in the studing of this generators is spectral analysis of their spectra and behaviour of its resolvents.  We consider that the results presented in the present paper are an important step in this direction because of the spectra of these generators coincides with the spectra of operator-valued function $L(\lambda)$. 

The paper is divided into five sections.  In the first section a brief introduction to the subject is carry out. Other works in the cases $\theta=0$, $\theta=1$ and  $\theta \in (0, 1)$ are  also mentioned. The main results on the correct solvability and spectral analysis of  integro-differential equation   \eqref{system 1.1} in the case $\theta \in [0, 1]$ are formulated in the second and the third
 paragraph, respectively.  The two theorems in the second  paragraph are related with correct solvability in the Sobolev space.  The first of them was proved in our work \cite{PORVVV}.  The second is proved in this present paper. The first result in the third paragraph is on the structure of the spectrum in the case when the kernel $K(t)$ belongs to the Sobolev space $W^{1}_{1}(\mathbb{R}_+)$, while the second result corresponds to the case when the kernel $K(t)$ belongs to the space $L_1(\mathbb{R}_+)$, but do not belong $W^{1}_{1}(\mathbb{R}_+)$. A complete analysis  on the structure of the spectrum of the operator-valued function $L(\lambda)$ in the case when $K(t)\in L_1(\mathbb{R}_+)$, but $K(t) \notin W^{1}_{1}(\mathbb{R}_+)$  is provided in the section \ref{representation-of-spectra-in-graphic-in-the-case-divergente-series}.  The proof of the main results is given in the fourth section. 
The fifth section contains the proof of auxiliary assertions and lemmas. 

\section{Correct solvability}
Let $H$ be a separable Hilbert space and let $A$ be a self-adjoint positive operator in~$H$ with compact inverse. We denote by $\{e_n\}_{n=1}^{\infty}$ the orthonormal basis formed by the eigenvectors of~$A$ corresponding to its eigenvalues $a_n$ such that $Ae_n=a_ne_n, n \in \mathbb{N}$.
 The eigenvalues  $a_n$ are arranged in increasing order
 and counted according to multiplicity; i. e, $$1\leq a_1 \leq a_2 \leq \cdots \leq a_n<\cdots,$$
 where $a_n \to +\infty$ as $n \to +\infty$.

Let us consider the system \eqref{system 1.1}--\eqref{system 1.2} on the semi-axis $\mathbb{R}_{+}=(0, +\infty)$. It is assumed that the vector-valued function $A^{2-\theta}f(t)$ belongs to $L_{2, \gamma}(\mathbb{R}_+, H)$ for some $\gamma \geq 0$, where
\begin{align*}
\|f\|_{L_{2, \gamma}(\mathbb{R}_{+}, H)}\equiv \left(\int_{0}^{\infty}\mathrm{e}^{-2\gamma t}\|f(t)\|^{2}_{H} dt\right)^{1/2},\hspace{0.5cm}\gamma \geq 0. 
\end{align*} The scalar function $K(t)$ admits the representation $$K(t)=\sum_{k=1}^{\infty}c_k\mathrm{e}^{-\gamma_k t},$$ where $c_k>0$, $\gamma_{k+1}>\gamma_{k}>0$, $k \in \mathbb{N}$, $\gamma_k \to +\infty$ $(k \to +\infty)$. Moreover, we assume that
\begin{align*}
a) \hspace{0.3cm}\sum_{k=1}^{\infty}\frac{c_k}{\gamma_k}<1, \hspace{1cm} b)\hspace{0.3cm} \sum_{k=1}^{\infty} c_k<+\infty.
\end{align*}
Note that if $a)$ holds, then $K \in L_1(\mathbb{R}_{+})$ and $\|K\|_{L_1(\mathbb{R}_{+})} <1$. If  $a)$ and $b)$ are both satisfied, then the kernel $K$ belongs to Sobolev space $W_{1}^{1}(\mathbb{R}_{+})$.

We denote by  $W_{2,\gamma}^{2}(\mathbb{R}_{+}, A^2)$  the Sobolev space consisting of vector-functions on the semi-axis $\mathbb{R}_{+}=(0, +\infty)$ with values in~$H$; this
space will be equipped with the norm
\begin{align*}
\|u\|_{W_{2,\gamma}^{2}(\mathbb{R}_{+}, A^2)} \equiv \left(\int_{0}^{\infty}\mathrm{e}^{-2\gamma t}\left(\|u^{(2)}(t)\|^{2}_{H}+\|A^2 u(t)\|^{2}_{H}\right) dt\right)^{1/2}, \hspace{0.3cm}\gamma \geq 0.
\end{align*}
For a complete description of the space $W_{2,\gamma}^{2}(\mathbb{R}_{+}, A^2)$ and some of its properties we refer to the monograph \cite[Chap. I]{LM}. For $\gamma=0$ we write $W_{2, 0}^{2}(\mathbb{R}_{+}, A^2)\equiv W^{2}_{2}(\mathbb{R}_{+}, A^2)$.  
\begin{definition} A vector-valued function $u$ is called a {\em strong solution} of system \eqref{system 1.1}--\eqref{system 1.2} if, for some $\gamma \geq 0$, $u \in W_{2,\gamma}^{2}(\mathbb{R}_{+}, A^2)$ satisfies the equation \eqref{system 1.1} almost everywhere on the semi-axis $\mathbb{R}_{+}$ and $u$~also satisfies  the initial condition \eqref{system 1.2}.
\end{definition} 
The following result  was proved in the work \cite{PORVVV}. 
\begin{theorem} \label{theorem-about-solvability} Suppose that, for all $\theta \in [0, 1]$ and  for some $\rho_0 \geq 0$, $A^{2-\theta} f(t)$ belongs to $L_{2,\rho_0}(\mathbb{R}_{+}, H)$. 
\begin{enumerate} [\ 1)]
\item \label{condition1-solvability} If conditions $a)$ and $b)$  both hold, and $\varphi_0 \in H_2$, $\varphi_1 \in H_1$ for all $\theta \in [0, 1]$ then there is a $\widetilde{\rho}>\rho_0$ such that for any $\gamma > \widetilde{\rho}$, the initial value problem \eqref{system 1.1}--\eqref{system 1.2}
 has a unique solution in the Sobolev space $W_{2,\gamma}^{2}(\mathbb{R}_{+}, A^2)$ and this solution satisfies the estimate
\begin{align}
\|u\|_{W_{2,\gamma}^{2}(\mathbb{R}_{+}, A^2)} \leq d\left(\|A^{2-\theta} f\|_{L_{2,\gamma}(\mathbb{R}_{+}, H)}+\|A^2\varphi_0\|_{H}+\|A\varphi_1\|_{H}\right),
\end{align} where the constant $d$ is independent of the vector-valued function $f$ and the vectors $\varphi_0$, $\varphi_1$.
\item  \label{condition2-solvability} If condition $a)$ is satisfied, but  condition $b)$  does not hold {\rm (}i.e., $K(t)\notin W_{1}^{1}(\mathbb{R}_{+}))$ and $\varphi_0 \in H_{2+\theta}$, $\varphi_1 \in H_{1+\theta}$ for all $\theta \in (0, 1]$ then there is a $\widetilde{\rho}>\rho_0$ such that for any $\gamma > \widetilde{\rho}$
the initial value problem \eqref{system 1.1}--\eqref{system 1.2} has a~unique solution in  Sobolev space $W_{2,\gamma}^{2}(\mathbb{R}_{+}, A^2)$ and this solution satisfies the estimate
\begin{align}
\|u\|_{W_{2,\gamma}^{2}(\mathbb{R}_{+}, A^2)} \leq d\left(\|A^{2-\theta} f\|_{L_{2,\gamma}(\mathbb{R}_{+}, H)}+\|A^{2+\theta}\varphi_0\|_{H}+\|A^{1+\theta} \varphi_1\|_{H}\right),
\end{align} where the constant $d$ is independent of the vector-valued function $f$ and the vectors $\varphi_0$, $\varphi_1$.
\end{enumerate}
\end{theorem}

We give a result on the correct solvability of the initial value problem  \eqref{system 1.1}--\eqref{system 1.2} in the Sobolev space $W_{2}^{2}((0, T), A^2)$ for any $T>0$. The space $W_{2}^{2}((0, T), A^2)$ is provided with the norm
\begin{align*}
\|u\|_{W_{2}^{2}((0, T), A^2)} \equiv \left(\int_{0}^{T}\left(\|u^{(2)}(t)\|^{2}_{H}+\|A^2 u(t)\|^{2}_{H}\right) dt\right)^{1/2}.
\end{align*}

\begin{theorem} \label{correct-solvability-theorem-abound-interval} Suppose that, for all $\theta \in [0, 1]$, vector-valued function $A^{2-\theta} f(t)$ belongs to space $L_{2}((0, T), H)$. Then
\begin{enumerate} [\ 1)]
\item If the conditions of part 1 of the theorem  \ref{theorem-about-solvability} are satisfied, then for an arbitrary $T>0$ the initial value problem \eqref{system 1.1}--\eqref{system 1.2} has a unique solution $u(t)$, which belongs to the Sobolev space $W_{2}^{2}((0, T), A^2)$, and for  that solution the following estimate holds
\begin{align*}
\|u\|_{W_{2}^{2}((0, T), A^2)} \leq d(T)\left(\|A^{2-\theta} f\|_{L_{2}((0, T), H)}+\|A^2\varphi_0\|_{H}+\|A\varphi_1\|_{H}\right),
\end{align*}  with a positive constant $d(T)$ that does not depend on the vector-function $f$ and vectors $\varphi_0$, $\varphi_1$.
\item  If the conditions of part 2 of the theorem \ref{theorem-about-solvability} are satisfied, then for an arbitrary $T>0$ the initial value problem \eqref{system 1.1}--\eqref{system 1.2} has a unique solution $u(t)$, which belongs to the Sobolev space $W_{2}^{2}((0, T), A^2)$, and for  that solution the following estimate holds 
\begin{align*}
\|u\|_{W_{2}^{2}((0, T), A^2)} \leq   d(T)\left(\|A^{2-\theta} f\|_{L_{2}((0, T), H)}+\|A^{2+\theta}\varphi_0\|_{H}+\|A^{1+\theta} \varphi_1\|_{H}\right),
\end{align*} with a positive constant $d(T)$ that does not depend on the vector-function $f$ and vectors $\varphi_0$, $\varphi_1$.
\end{enumerate}
\end{theorem}

It is relevant to note that several results about correct solvability in the another functional spaces on finite interval $(0, T)$, $T>0$, and under another understanding of solutions of the initial value problem \eqref{system 1.1}--\eqref{system 1.2} were obtained in   \cite{GMJ}  (see also \cite{MRWD,P}). 
\section{Spectral analysis}
Consider the operator-valued function
\begin{align}
L(\lambda)= \lambda^2I+A^2-\widehat{K}(\lambda)A^{2\theta}, \hspace{1cm} 0\leq \theta \leq 1,
\end{align} which is the symbol of integrodifferential equation \eqref{system 1.1}, where
$$\widehat{K}(\lambda)=\sum_{k=1}^{\infty}\frac{c_k}{\lambda+\gamma_k}$$ is the Laplace transform of the kernel $K(t)$, the operator $I$ is the unit operator acting in a separable Hilbert space $H$. 

Let us consider the restriction of the operator-valued function $L(\lambda)$ on a one-dimensional subspace spanned by a vector $e_n$:
\begin{align} \label{meromorfonaya-funktsia}
\ell_n(\lambda):=(L(\lambda)e_n, e_n)= \lambda^2+ a^{2}_{n} \left(1- \frac{1}{a^{2(1-\theta)}_n}\sum_{k=1}^{\infty} \frac{c_k}{\lambda+\gamma_k}\right), \hspace{0.2cm}  0\leq \theta \leq 1.
\end{align} 

In what follows it is assumed that the following condition is satisfied
\begin{align}\label{suprem-de-las-gammas}
\sup_{k \in \mathbb{N}}\{\gamma_k(\gamma_{k+1}-\gamma_k)\}=+\infty. 
\end{align} The condition \eqref{suprem-de-las-gammas} means that the elements of the sequence $\{\gamma_k\}^{\infty}_{k=1}$, can not approach each other too quickly.

\subsection{Theorem on the location of spectrum of operator-valued functions}
\begin{definition} By the resolvent set $R(\lambda)$ of operator-valued function $L(\lambda)$ it is understood as the set of all values $\lambda \in \mathbb{C}$, for which the operator-valued function $L^{-1}(\lambda)$ exists and it is bounded. The complement of set $R(\lambda)$ in the complex plane, i. e., $\sigma(L)=\{\mathbb{C}\setminus R(\lambda)\}$, is called the spectrum of the operator-valued function $L(\lambda)$.
\end{definition}
\begin{theorem}  \label{spectr-operator-function-lying-on-the-left-right} Suppose  that $\sum_{k=1}^{\infty}\frac{c_k}{\gamma_k}<a^{2(1-\theta)}_1$  is satisfied. Then the spectrum of operator-valued function $L(\lambda)$ is contained in the left half-plane $\{\lambda \in \mathbb{C}: \re \lambda<0\}$. 
\end{theorem}  

\begin{remark}
If $\sum_{k=1}^{\infty}\frac{c_k}{\gamma_k}>a^{2(1-\theta)}_n$, $n=1, \dots, N_0$  is satisfied, then  on the right half-plane  lies $N_0$ positive eigenvalues of the operator-valued function  $L(\lambda)$. 
\end{remark}

\subsection{Theorem on the structure of spectra in the case when the kernel $K(t)$ belongs to the Sobolev space $W_{1}^{1}(\mathbb{R}_+)$} \label{analysis-about-spectrum-in-the-easy-case}
\begin{theorem} \label{theorem-about-spectral-anal-1} Suppose that the following conditions  \eqref{suprem-de-las-gammas},  $a)$, $b)$ and $a_1 \geq 1$ hold. Then, for each fixed $n \in \mathbb{N}$, the set of zeros of meromorphic function $\ell_n(\lambda)$ is the union of a countable set of real zeros $\{\lambda_{n, k}(\theta)|n, k \in \mathbb{N}\}$ satisfying the inequalities
\begin{align}\label{desigualdad-ceros-finitos}
\cdots-\gamma_k<\lambda_{n, k}(\theta) < \cdots <-\gamma_{1}<\lambda_{n, 1}(\theta)<0,  \hspace{0.3cm}  \mbox{and}  \hspace{0.3cm} \lim_{n\to +\infty}\lambda_{n, k}(\theta)=-\gamma_k,
\end{align} and also pairs of zeros $\lambda^{\pm}_{n}(\theta)$, which, for a sufficiently large $n \in \mathbb{N}$, are non-real, complex-conjugate $\lambda^{+}_{n}(\theta)= \overline{\lambda^{-}_{n}}(\theta)$ and  asymptotically represented in the form
 \begin{align}\label{asymptota-sum-finita-para-mu1}
\lambda^{\pm}_{n}(\theta)=-\frac{1}{2}\frac{1}{a_{n}^{2(1-\theta)}}\sum \limits_{k=1}^{\infty}c_k+O\left(\frac{1}{a^{4-2\theta}_{n}}\right)\pm i\left(a_{n}+O\left(\frac{1}{a^{3-2\theta}_{n}}\right)\right), \hspace{0.1cm} a_n \to +\infty. 
\end{align}
\end{theorem}  
\begin{remark} In the relation \eqref{asymptota-sum-finita-para-mu1}, the subordinate summands containing symbols $O\left(\frac{1}{a^{k}_{n}}\right)$ are written separately for the real and imaginary parts of zeros $\lambda^{\pm}_{n}(\theta)$.
\end{remark}
\begin{corollary} If the conditions of Theorem \ref{theorem-about-spectral-anal-1} are satisfied, the spectrum $\sigma (L)$ of the operator-valued function  $L(\lambda)$ coincides with the closure of set of zeros $\{\lambda^{\pm}_{n}(\theta)\}^{\infty}_{n=1}$ and $\{\lambda_{n, k}(\theta)\}^{\infty, \infty}_{n, k=1}$ of meromorphic function $\ell_n(\lambda)$, i. e., the spectrum $\sigma (L)$ is represented by
\begin{align*}
 \sigma(L)=\overline{\left(\bigcup^{\infty}_{n=1}\bigcup^{\infty}_{k=1}\lambda_{n, k}(\theta)\right) \cup \left(\bigcup^{\infty}_{n=1} \lambda^{\pm}_{n}(\theta)\right)}.
\end{align*}
\end{corollary} 
 It is observed that if the conditions of theorem \ref{theorem-about-spectral-anal-1} are satisfied, then in the case when $\theta \in [0, 1)$ the non-real parts of complex-conjugate roots  $\lambda^{\pm}_{n}(\theta)$, when  $n \to +\infty$, asymptotically approach the imaginary axis (see figure 1).  In the case $\theta=1$, the non-real parts of complex conjugate roots $\lambda^{\pm}_{n}(\theta)$, when  $n \to +\infty$,  asymptotically approach a line parallel to the imaginary axis (for more details, see the works \cite{VR, VR1, VR2} and chapter 3 of  monograph \cite{vlasovmrau2016}). Thus, when $\theta \in [0, 1)$, the non-real spectrum of the operator-valued function $L(\lambda)$ is close to the spectrum of the abstract wave equation (when $K(t)\equiv 0$). 
 
 The structure of the real spectrum $\sigma_{\mathbb{R}}$ when $\theta \in [0, 1)$ also differs from the case $\theta=1$, since when $\theta=1$ the real zeros $\lambda_{n, k}(\theta)$ tend to real zeros $x_k$ of function $$f(\lambda)=1-\sum_{k=1}^{\infty}\frac{c_k}{\lambda+\gamma_k}$$ (see  \cite{VR, VRShamaev, VR2} for more details). 

\setlength{\unitlength}{1cm}
\begin{picture}(1, 1)(-1, 3)     
\put(2, 0){\vector(1,0){8}}
\put(10.9,0){\makebox(0,0){$\re \lambda$}}
 \put(8,-3){\vector(0,1){6}}
\put(8.9,3){\makebox(0,0){$\im \lambda$}}
\put(6.4,-3.8){\makebox(0,0){{\bf Figure 1}: Structure of spectrum of operator-valued  function $L(\lambda)$}}
\put(5.3,-4.3){\makebox(0,0){ when the kernel $K(t) \in W_{1}^{1}(\mathbb{R}_+)$.}}
\put(3,-1.5){\makebox(0,0){$\lambda_{n, k}(\theta)$}}
\put(3, -1.3){\vector(1,1){1}}
\put(2,0.6){\makebox(0,0){$-\infty$}}
\put(2.4, 0.3){\vector(-1,0){0.5}}
\put(2,0){{\circle*{0.15}}}
\put(2.1,0){{\circle*{0.15}}}
\put(2.2,0){{\circle*{0.15}}}
\put(2.3,0){{\circle*{0.15}}}
\put(2.4,0){{\circle*{0.15}}}
\put(2.8,0){{\circle*{0.15}}}
\put(2.9,0){{\circle*{0.15}}}
\put(3,0){{\circle*{0.15}}}
\put(3.1,0){{\circle*{0.15}}}
\put(3.2,0){{\circle*{0.15}}}
\put(3.8,0){{\circle*{0.15}}}
\put(3.9,0){{\circle*{0.15}}}
\put(4,0){{\circle*{0.15}}}
\put(4.1,0){{\circle*{0.15}}}
\put(4.2,0){{\circle*{0.15}}}
\put(4.8,0){{\circle*{0.15}}}
\put(4.9,0){{\circle*{0.15}}}
\put(5,0){{\circle*{0.15}}}
\put(5.1,0){{\circle*{0.15}}}
\put(5.2,0){{\circle*{0.15}}}
\put(5.8,0){{\circle*{0.15}}}
\put(5.9,0){{\circle*{0.15}}}
\put(6,0){{\circle*{0.15}}}
\put(6.1,0){{\circle*{0.15}}}
\put(6.2,0){{\circle*{0.15}}}
\put(4.9,2.5){\makebox(0,0){$\lambda^{\pm}_{n}(\theta)$}}
\put(5,2){\vector(1,-1){1}}
\put(7.3,3){{\circle*{0.15}}}
\put(7.1,2.2){{\circle*{0.15}}}
\put(6.7,1.4){{\circle*{0.15}}}
\put(6.2,0.8){{\circle*{0.15}}}
\put(5.6,0.4){{\circle*{0.15}}}
\put(4.9,0.15){{\circle*{0.15}}}
\put(4.9,-0.15){{\circle*{0.15}}}
\put(5.6,-0.4){{\circle*{0.15}}}
\put(6.2,-0.8){{\circle*{0.15}}}
\put(6.7,-1.4){{\circle*{0.15}}}
\put(7.1,-2.2){{\circle*{0.15}}}
\put(7.3,-3){{\circle*{0.15}}}
\put(8.5,1.0){\makebox(0,0){$ia_n$}}
\put(8.5,-1.0){\makebox(0,0){$-ia_n$}}
\end{picture}

\vspace{8cm}

\begin{remark} In this situation, the solution of homogeneous integrodifferential equation \eqref{system 1.1} can not decay exponentially. Indeed, if the solutions decayed exponentially like $\mathrm{e}^{-\alpha t}$ we would have a vertical strip $\{\lambda:\alpha<\re \lambda<0\}$ free from spectra of operator-valued function $L(\lambda)$.  
\end{remark}
\subsection{Theorem on the structure of the spectrum in the case when the kernel  $K(t)$  belongs to space   $L_1(\mathbb{R}_+)$, but does not belong to space $W_{1}^{1}(\mathbb{R}_+)$} 

{\bf Condition 1}. Suppose that the sequences  
 $\{c_{k}\}^{\infty}_{k=1}$ and $\{\gamma_{k}\}^{\infty}_{k=1}$ have the following representation
 \begin{align*}
 c_k&=\frac{\mathcal{A}}{k^{\alpha}} + O\left( \frac{1}{k^{\alpha +1}}\right), \hspace{0.2cm}k \in \mathbb{N}\\  
 \gamma_k&=\mathcal{B} k^{\beta} + O( k^{\beta-1}),  \hspace{0.2cm}k \in \mathbb{N},
\end{align*} where the constants $\mathcal{A}>0, \mathcal{B}>0$, $0<\alpha \leq 1$, $\alpha +\beta>1$, and, when $k \to +\infty$, the sequences $\{c_{k}\}^{\infty}_{k=1}$ and $\{\gamma_{k}\}^{\infty}_{k=1}$ satisfy the condition $\sum_{k=1}^{\infty}\frac{c_k}{\gamma_k}<1$. 
\begin{remark} If the following relation holds
\begin{align}
\gamma_{k+1}-\gamma_{k} \thickapprox k^{\beta-1}, \hspace{0.5cm}k \to +\infty,
\end{align}  then for $\beta>1/2$ the condition \eqref{suprem-de-las-gammas} is satisfied. 
\end{remark}
\begin{remark} It is noticed that if the {\bf Condition 1} is fulfilled, the kernel $K(t)$ have a singularity at $t=0$,  because  $$K(0)=\sum_{k=1}^{\infty}c_k=+\infty.$$ 

The following theorem represents the asymptotics of a pair of complex-conjugate zeros $\lambda^{\pm}_{n}$,  $\lambda^{+}_{n}= \overline{\lambda^{-}_{n}}$ in the case when the condition  $b)$ is not satisfied.
\end{remark}
\begin{theorem} \label{theorem-when-condition-not-hold} Suppose that $\beta>1/2$, $a_{1}\geq 1$, {\bf Condition 1}  and  $a)$ are satisfied. Then, for each fixed $n\in \mathbb N$, the set of zeros of a meromorphic function $\ell_n(\lambda)$ is the union of a countable set of real zeros satisfying the inequalities \eqref{desigualdad-ceros-finitos} with a pairs of zeros $\lambda^{\pm}_{n}(\theta, r)$, which, for a sufficiently large $n \in \mathbb{N}$, are non-real, complex-conjugate  $\lambda^{+}_{n}(\theta, r)= \overline{\lambda^{-}_{n}}(\theta, r)$ and asymptotically represented, when $a_n \to +\infty$, in the following form
  \begin{align}
  \lambda^{\pm}_{n}(\theta, r)&=-\frac{\mathcal{A}D_1 \mathcal{B}^{r-1}}{\beta a^{n_1(\theta, r)}_n}\pm i\left(a_n+\frac{\mathcal{A}D_2\mathcal{B}^{r-1} }{\beta a^{n_1(\theta, r)}_n}\right)+O\left(\frac{1}{a^{n_2(\theta, r)}_n}\right), \mbox{$r\in (0, \frac{1}{2}) \wedge \theta \in [\frac{1}{2}, 1)$},\label{asimptotas-de-ceros-infinitos-1} \\
  \lambda^{\pm}_{n}(\theta, r)&=-\frac{\mathcal{A}D_1\mathcal{B}^{r-1}}{\beta a^{n_1(\theta, r)}_n}\pm i\left(a_n+\frac{\mathcal{A}D_2\mathcal{B}^{r-1}}{\beta a^{n_1(\theta, r)}_n}\right)+O\left(\frac{1}{a^{2(1-\theta)}_n}\right), \mbox{$r\in [\frac{1}{2}, 1) \vee \theta \in (0, \frac{1}{2})$},\label{asimptotas-de-ceros-infinitos-2} \\
 \lambda^{\pm}_{n}(\theta, r)&=-\frac{1}{2}\frac{\mathcal{A}}{\beta}\frac{\ln a_n}{a^{2(1-\theta)}_n}\pm ia_n+O\left(\frac{1}{a^{2(1-\theta)}_n}\right), \hspace{5cm}\mbox{$r=1$},\label{asimptotas-de-ceros-infinitos-3} 
\end{align} where $n_1(\theta, r):=r+2\left(\frac{1}{2}-\theta\right)$, $n_2(\theta, r):=\min \{2(1-\theta), 2r+3-4\theta\}$, the parameter $r:=\frac{\alpha+\beta-1}{\beta}$,  $\alpha$ and $\beta$ are such that $\alpha \in (0, 1]$, $\alpha +\beta>1$, the constants  $\mathcal{A}>0, \mathcal{B}>0$ and constants $D_1$, $D_2$ are defined as follows
 \begin{align*} 
D_1:=\frac{\pi}{2}\frac{\cos\left(\frac{\pi}{2} (r+1)\right)}{\sin(\pi r)}, \hspace{1cm} D_2:=-\frac{\pi}{2} \frac{\sin\left(\frac{\pi}{2} (r+1)\right)}{\sin(\pi r)}. 
\end{align*}
\end{theorem}
 \begin{corollary} If the conditions of Theorem \ref{theorem-when-condition-not-hold} are satisfied, the spectrum $\sigma (L)$ of the operator-valued function $L(\lambda)$ coincides with the closure of the set of zeros $\{\lambda^{\pm}_{n}\}^{\infty}_{n=1}(\theta, r)$
and $\{\lambda_{n, k}\}^{\infty, \infty}_{n, k=1}(\theta)$ of meromorphic function $\ell_n(\lambda)$, i. e., the spectrum $\sigma (L)$ is  represented by 
\begin{align*}
 \sigma(L)=\overline{\left(\bigcup^{\infty}_{n=1}\bigcup^{\infty}_{k=1}\lambda_{n, k}(\theta)\right) \cup \left(\bigcup^{\infty}_{n=1} \lambda^{\pm}_{n}(\theta, r)\right)}.
\end{align*}  
\end{corollary}
\begin{remark} The case $\theta=1$ was studied in detail in \cite{VR, VR1,VR2}, as well as in monograph \cite[chap. 3]{vlasovmrau2016}. For $\theta=1$, the proposed theorem  \ref{theorem-when-condition-not-hold} becomes theorem 3 of the paper \cite{VR}.
\end{remark}

\section{Proof of the main results} \label{section-proof-about-spectr-operator-function}
\subsection{Proof of theorem \ref{correct-solvability-theorem-abound-interval}}

We carry out the proof of the mentioned theorem at first for zero initial conditions $u(+0)=\varphi_0=0$, $u^{(1)}(+0)=\varphi_1=0$. We extend the function $f(t)$, given on the segment $[0, T]$,  as zero function on the semiaxis $(T, +\infty)$. Let us denote the extension of the function $f(t)$ by $\widetilde{f}(t)$:  
 
\[ \widetilde{f}(t) = \left\{ \begin{array}{ll}
        f(t) & \mbox{if $t \in [0, T]$},\\
       0 & \mbox{if $t>T$}.\end{array} \right. \]
Then the function   $\widetilde{f}(t)$   will satisfy the conditions of the theorem  \ref{theorem-about-solvability} with the vector-valued function $A^{2-\theta}\widetilde{f}(t) \in L_{2, \gamma}(\mathbb{R}_+, H)$, $\gamma>0$, and the inequality holds
\begin{align} \label{inequality-when-L-2-2-gamma-less-L-2}
\left\|A^{2-\theta}\widetilde{f}(t) \right\|_{L_{2, \gamma}(\mathbb{R}_+, H)}\leq \left\|A^{2-\theta}f(t) \right\|_{L_{2}((0, T), H)}. 
\end{align} Further, according to the theorem \ref{theorem-about-solvability}, for the function $\widetilde{f} (t) $ there exists a unique solution $u(t)$, which belongs to the Sobolev space $W_{2, \gamma}^{2} (\mathbb{R}_+, A^2) $, for $\gamma>\tilde {\rho}> 0$, and satisfying the estimate 
\begin{align}\label{inequality-when-W-2-2-gamma-less-W-2-2}
\left\|u \right\|_{W_{2, \gamma}^{2}(\mathbb{R}_+, A^2)} \leq d_0 \left\|A^{2-\theta}\widetilde{f}(t) \right\|_{_{L_{2, \gamma}(\mathbb{R}_+, H)}}, \gamma>0. 
\end{align} It is not difficult to note that the following chain of inequalities holds: 
\begin{align} \label{chain-of-inequalities-when-w-2-2-and-w-2-2-gamma}
 \|u\|^2_{W^{2}_{2, \gamma}(\mathbb{R}_+, A^2)}&=  \int_{0}^{\infty}\mathrm{e}^{-2\gamma t}\left(\|u^{(2)}(t)\|^2_H+\|A^2u(t)\|^2_H \right)dt \nonumber \\ 
&\geq \int_{0}^{T}\mathrm{e}^{-2\gamma t}\left(\|u^{(2)}(t)\|^2_H+\|A^2u(t)\|^2_H \right)dt\nonumber \\ 
&\geq \mathrm{e}^{-2\gamma T}\int_{0}^{T}\left(\|u^{(2)}(t)\|^2_H+\|A^2u(t)\|^2_H \right)dt\nonumber\\
&= \mathrm{e}^{-2\gamma T}\|u\|^2_{W^{2}_{2}((0, T), A^2)}. 
\end{align} Combining the inequalities  \eqref{inequality-when-L-2-2-gamma-less-L-2},  \eqref{inequality-when-W-2-2-gamma-less-W-2-2} and \eqref{chain-of-inequalities-when-w-2-2-and-w-2-2-gamma} we arrive at the desired inequality
\begin{align} \label{desired-inequality-for-space-without-weigth}
\|u\|_{W^{2}_{2}((0, T), A^2)}&=\left(\int_{0}^{T}\left(\|u^{(2)}(t)\|^2_H+\|A^2u(t)\|^2_H \right)dt\right)^{1/2}\nonumber\\
 &\leq d_0\mathrm{e}^{\gamma T}\|A^{2-\theta}f\|_{L_{2}((0, T), H)}. 
\end{align} Let us consider the case of non-zero initial conditions $\varphi_0$ and  $\varphi_1$. As in the proof of the theorem \ref{theorem-about-solvability} we seek the solution of the initial value problem \eqref{system 1.1}--\eqref{system 1.2} in the form 
\begin{align}\label{solutions-in-the-form-without-weigth-gamma}
u(t)=\cos(At)\varphi_0+A^{-1}\sin (At)\varphi_1+\omega(t). 
\end{align} Then for the function  $\omega(t)$ we obtain the following initial value  problem: 
\begin{align}
\frac{d^2}{dt^2}\omega(t)+&A^2 \omega(t)- \int_{0}^{t}K(t-s)A^{2\theta}\omega(s)ds=f_1(t), \label{equation-w-when-0-T}\\
&\omega(+0)=\omega^{(1)}(+0)=0, \label{equation-w-when-0-T-inicial-conditions}
\end{align} where $f_1(t)=f(t)+h(t)$ and 
\begin{align*}
h(t)&=\int_{0}^{t}K(t-s)A^{2\theta}\left(\cos (As)\varphi_0+A^{-1}\sin (As)\varphi_1\right)ds, \\
A^{2-\theta}h(t)&=\int_{0}^{t}K(t-s)\left(A^{2+\theta}\cos (As)\varphi_0+A^{1+\theta}\sin (As)\varphi_1\right)ds. 
\end{align*} First, this assertion will be proven when the condition $1)$ of theorem \ref{theorem-about-solvability} is satisfied, i. e., in the case when the series  $\sum \limits_{k=1}^{\infty} c_k$ is convergent. By the linearity of the equation \eqref{equation-w-when-0-T} the solution $\omega(t)$ can be represented in the form
\begin{align} \label{lineal-representation-of-omega-in-omega-1-and-omega-2}
\omega(t)= \omega_1(t)+\omega_2(t), 
\end{align} where $\omega_1(t)$ is the solution of the problem corresponding to the right-hand side of $h(t)$, while $\omega_2(t)$ is the solution of the problem corresponding to the right-hand side of  $f(t)$. The estimate of function $\omega_2(t)=u(t)$ already was obtained in
\eqref{desired-inequality-for-space-without-weigth}. We now estimate the function $\omega_1(t)$. We note that in the proof of the theorem \ref{theorem-about-solvability} in the case when the condition  $1)$ is satisfied the following inequality was established
\begin{align} \label{estimate-on-the-space-L-and-gamma}
\|A^{2-\theta}h\|_{L_{2, \gamma}(\mathbb{R}_+, H)}\lesssim \left(\sum_{k=1}^{\infty}c_k\right)\left(\|A^{1+\theta}\varphi_0\|+\|A^{\theta}\varphi_1\|\right). 
\end{align} From the inequality \eqref{estimate-on-the-space-L-and-gamma} and according to the proof of the theorem \ref{theorem-about-solvability} the following estimate is valid
 \begin{align}\label{estimate-under-omega-1-in-w-2-2-gamma-0}
\|\omega_1\|_{W_{2, \gamma}^{2}(\mathbb{R}_+, A^2)}\leq d_1\|A^{2-\theta}h\|_{L_{2, \gamma}(\mathbb{R}_+, H)}\leq d_2 \left(\|A^{1+\theta}\varphi_0\|+\|A^{\theta}\varphi_1\|\right),
\end{align}  with constants $d_1$ and $d_2$ not depending on the vectors $\varphi_0$ and $\varphi_1$.  In turn, according to what was proved before (see inequality  \eqref{chain-of-inequalities-when-w-2-2-and-w-2-2-gamma}) for function $\omega_1(t)$ the estimate is satisfied 
\begin{align}\label{estimate-under-omega-1-in-w-2-2-gamma-1}
\|\omega_1\|_{W_{2}^{2}((0, T), A^2)}\leq \mathrm{e}^{\gamma T}\|\omega_1\|_{W_{2, \gamma}^{2}(\mathbb{R}_+, A^2)}. 
\end{align} Therefore, from the estimates \eqref{estimate-under-omega-1-in-w-2-2-gamma-0} and \eqref{estimate-under-omega-1-in-w-2-2-gamma-1} we obtain the required inequality
\begin{align}\label{estimate-under-omega-1-in-w-2-2-gamma-2}
\|\omega_1\|_{W_{2}^{2}((0, T), A^2)}\leq d_1(T) \left(\|A^{1+\theta}\varphi_0\|+\|A^{\theta}\varphi_1\|\right). 
\end{align}  Consequently, from the representation \eqref{lineal-representation-of-omega-in-omega-1-and-omega-2} and inequalities  \eqref{desired-inequality-for-space-without-weigth}, \eqref{estimate-under-omega-1-in-w-2-2-gamma-2} with $u(t)=\omega_2(t)$ we obtain the estimate
\begin{align}\label{estimate-under-omega-1-in-w-2-2-gamma-3}
\|\omega_1\|_{W_{2}^{2}((0, T), A^2)}\leq d_2(T) \left(\|A^{2-\theta}f\|_{L_{2}((0, T), H)}+\|A^{1+\theta}\varphi_0\|+\|A^{\theta}\varphi_1\|\right). 
\end{align} To complete the proof, we give an estimate of the vector function
\begin{align*}
v(t)=\cos (At)\varphi_0+A^{-1}\sin (At)\varphi_1
\end{align*} in the Sobolev space $W_{2}^{2}((0, T), A^2)$, for any $T>0$.  

 According to known inequalities $\|\cos(At)\|_H\leq 1$, $\|\sin(At)\|_H\leq1$ for $\varphi_0 \in H_2= \dom (A^2)$ and $\varphi_1 \in H_1= \dom(A)$ the following inequalities are obtained
 \begin{align}
\|\cos(At)\varphi_0\|^{2}_{W_{2}^{2}((0, T), A^2)} &= \int \limits_{0}^{T}\left(\|A^2\cos(At)\varphi_0\|^{2}_{H}+\|A^2 \cos(At)\varphi_0\|^{2}_{H}\right) dt\nonumber\\
&\leq 2\int \limits_{0}^{T}\|A^2\varphi_0\|^{2}_{H}dt=2T\|A^2\varphi_0\|^{2}_{H},\label{inequality-cosA-1}\\
\|A^{-1}\sin(At)\varphi_1\|^{2}_{W_{2}^{2}((0, T), A^2)} &= \int\limits_{0}^{T}\left(\|A\sin(At)\varphi_1\|^{2}_{H}+\|A \sin(At)\varphi_1\|^{2}_{H}\right) dt\nonumber\\
&\leq 2\int \limits_{0}^{T}\|A\varphi_1\|^{2}_{H} dt=2T\|A\varphi_1\|^{2}_{H}. \label{inequality-sinA-1}
\end{align} From the inequalities \eqref{inequality-cosA-1} and \eqref{inequality-sinA-1}  the following estimate follows 
\begin{align} \label{estimate-under-v-with-initial-conditions}
\|v\|_{W^2_{2} ((0, T), A^2)}\leq \sqrt{2T}\left(\|A^{2}\varphi_0\|_{H}+\|A\varphi_1\|_{H}\right). 
\end{align} On the basis of inequalities \eqref{estimate-under-omega-1-in-w-2-2-gamma-3}, \eqref{estimate-under-v-with-initial-conditions} and representation \eqref{solutions-in-the-form-without-weigth-gamma} the required estimate is obtained 
\begin{align*} 
\|u(t)\|_{W^2_{2} ((0, T), A^2)}\leq D(T)\left(\|A^{2-\theta}f\|_{L_{2} ((0, T), H)}+\|A^{2}\varphi_0\|_{H}+\|A\varphi_1\|_{H}\right),
\end{align*}  where the positive constant $D(T)$ is  independent of the vector function $f$ and vectors $\varphi_0$, $\varphi_1$. 

Consider the case when the series $\sum_{k=1}^{\infty}c_k$ is divergent, i. e., when the condition  $b)$
is not satisfied, but $\sum\limits_{k=1}^{\infty}\frac{c_k}{\gamma_k}<1$ is satisfied. Repeating the reasoning of the part  $1)$ of theorem \ref{correct-solvability-theorem-abound-interval}   we can estimate the function $\omega_2(t)$: 
\begin{align}\label{omega-2-when-the-series-under-c-k-diverges}
\|\omega_2\|_{W^2_{2} ((0, T), A^2)}\leq d_0 \mathrm{e}^{\gamma T} \left(\|A^{2-\theta}f\|_{L_{2} ((0, T), H)}\right). 
\end{align} For the vector-function $v(t)$ the estimate \eqref{estimate-under-v-with-initial-conditions} is valid.  To obtain the estimate of function $w_1$ it is necessary the proof of  the following lemma: 
\begin{lemma} \label{lemma-when-series-under-c-k-diverges}  Under the assumptions that  the series  $\sum_{k=1}^{\infty}\frac{c_k}{\gamma_k}$ is convergent, the following estimate is true 
\begin{align} \label{estimate-when-series-under-c-k-diverges-and-series-under-c-k-gamma-converge}
\|A^{2-\theta}h\|_{L_{2, \gamma} (\mathbb{R}_+, H)}\leq \operatorname{const} \left(\sum \limits_{k=1}^{\infty}\frac{c_k}{\gamma_k}\right) \left(\|A^{2+\theta}\varphi_0\|_{H}+\|A^{1+\theta}\varphi_1\|_{H}\right). 
\end{align}
\end{lemma}

\textbf{Proof.}  The following chain of inequalities is quickly verified 
 \begin{align} \label{estimate-operator-A-2-theta-when-series-under-c-k-diverges}
\|A^{2-\theta}h&\|_{\mathcal{L}_{2,\gamma}}\leq\left\|\sum \limits_{k=1}^{\infty}s_k A^{2+\theta}\left(-\gamma_k\delta_k+A\sin(At)\right)\varphi_0\right\|_{\mathcal{L}_{2,\gamma}}+\left\|\sum \limits_{k=1}^{\infty} s_k A^{1+\theta}\left(A\delta_k+\gamma_k\sin(At)\right)\varphi_1\right\|_{\mathcal{L}_{2,\gamma}}\nonumber\\
&\leq\left\|\cos (At) \sum \limits_{k=1}^{\infty} \frac{c_k}{\gamma_k} \gamma^2_k[A^2+\gamma^2_kI]^{-1} A^{2+\theta}\varphi_0 \right\|_{\mathcal{L}_{2,\gamma}}+ \left\|\sum \limits_{k=1}^{\infty} \frac{c_k }{\gamma_k}\frac{\gamma^2_k}{\mathrm{e}^{\gamma_kt}} [A^2+\gamma^2_kI]^{-1} A^{2+\theta}\varphi_0 \right\|_{\mathcal{L}_{2,\gamma}}\nonumber\\
&+\left\|\sin(At)\sum \limits_{k=1}^{\infty} \frac{c_k}{\gamma_k} \gamma_k A [A^2+\gamma^2_kI]^{-1} A^{2+\theta}\varphi_0 \right\|_{\mathcal{L}_{2,\gamma}} +\left\|\sum \limits_{k=1}^{\infty} \frac{c_k }{\gamma_k} \frac{\gamma_k}{\mathrm{e}^{\gamma_kt}} [A^2+\gamma^2_kI]^{-1} A^{2+\theta}\varphi_1 \right\|_{\mathcal{L}_{2,\gamma}}\nonumber\\
& +\left\|\cos (At)\sum \limits_{k=1}^{\infty}  \frac{c_k}{\gamma_k}\gamma_k [A^2+\gamma^2_kI]^{-1}  A^{2+\theta}\varphi_1 \right\|_{\mathcal{L}_{2,\gamma}}+\left\|\sin(At)\sum \limits_{k=1}^{\infty}  \frac{c_k }{\gamma_k}\gamma^2_k[A^2+\gamma^2_kI]^{-1}  A^{1+\theta}\varphi_1 \right\|_{\mathcal{L}_{2,\gamma}}\nonumber\\
&\leq\left\|\sum \limits_{k=1}^{\infty}   \frac{c_k }{\gamma_k}\gamma^2_k[A^2+\gamma^2_kI]^{-1}  A^{2+\theta}\varphi_0 \right\|_{\mathcal{L}_{2,\gamma}}+  \left\|\sum \limits_{k=1}^{\infty}\frac{c_k}{\gamma_k} \frac{\gamma^2_k}{\mathrm{e}^{\gamma_kt}} [A^2+\gamma^2_kI]^{-1}   A^{2+\theta}\varphi_0 \right\|_{\mathcal{L}_{2,\gamma}}\nonumber\\
&+\left\|\sum \limits_{k=1}^{\infty}   \frac{c_k }{\gamma_k}\gamma_kA [A^2+\gamma^2_kI]^{-1} A^{2+\theta}\varphi_0 \right\|_{\mathcal{L}_{2,\gamma}}+\left\|\sum \limits_{k=1}^{\infty}  \frac{c_k}{\gamma_k}  \frac{\gamma_k}{\mathrm{e}^{\gamma_kt}} A [A^2+\gamma^2_kI]^{-1} A^{1+\theta}\varphi_1 \right\|_{\mathcal{L}_{2,\gamma}}\nonumber\\
&+\left\|\sum \limits_{k=1}^{\infty}  \frac{c_k }{\gamma_k}\gamma_k A[A^2+\gamma^2_kI]^{-1}  A^{1+\theta}\varphi_1 \right\|_{\mathcal{L}_{2,\gamma}}+\left\|\sum \limits_{k=1}^{\infty} \frac{c_k}{\gamma_k} \gamma^2_k[A^2+\gamma^2_kI]^{-1}  A^{1+\theta}\varphi_1 \right\|_{\mathcal{L}_{2,\gamma}}
\end{align} where $\delta_k=\mathrm{e}^{-\gamma_kt}I-\cos (At)$, $s_k=c_k(A^2+\gamma^2_kI)^{-1}$ and  
$\mathcal{L}_{2,\gamma}:=L_{2,\gamma} (\mathbb{R}_+, H)$.

 To prove the lemma \ref{lemma-when-series-under-c-k-diverges} we need the following easily verifiable proposition. 
\begin{proposition} \label{proposition-about-three-termins-when-auto-adjoint-operator} Under the assumptions established on the  sequences $\{\gamma_k\}_{k=1}^{\infty}$, $\gamma_k \to +\infty$ $(k\to +\infty)$ and on the operator $A$ the following inequalities hold
\begin{align}
\|\gamma^2_k(A^2+\gamma^2_kI)^{-1}\|&< 1, \hspace{5mm} k\in \mathbb{N}, \label{norm-of-operator-A-and-gamma-k-1}\\
\|\gamma_k A(A^2+\gamma^2_kI)^{-1}\|&\leq \frac{1}{2}, \hspace{5mm} k\in \mathbb{N}, \label{norm-of-operator-A-and-gamma-k-2}
\end{align} 
\end{proposition}
The inequalities  \eqref{norm-of-operator-A-and-gamma-k-1} and \eqref{norm-of-operator-A-and-gamma-k-2} follow on from the inequalities
\begin{align*}
\frac{\gamma^2_k }{\gamma^2_k + a^2_n}&< 1, \hspace{5mm} n, k\in \mathbb{N}, \\
 \frac{\gamma_k a_n}{\gamma^2_k + a^2_n}&\leq  \frac{1}{2}, \hspace{5mm} n, k\in \mathbb{N}, 
\end{align*} and from the representation of operator $A$: 
$$A=\sum_{n=1}^{\infty}a_n (\cdot, e_n) e_n,$$ 
where $\{e_n\}_{n=1}^{\infty}$ is the orthonormal basis of eigenvectors of the operator $A$: $Ae_n=a_ne_n$. 

By the proposition \ref{proposition-about-three-termins-when-auto-adjoint-operator} and  the convergence of series 
$\sum_{k=1}^{\infty}\frac{c_k}{\gamma_k}$  each of the six terms in \eqref{estimate-operator-A-2-theta-when-series-under-c-k-diverges} can be estimated. Consequently, the estimate  \eqref{estimate-when-series-under-c-k-diverges-and-series-under-c-k-gamma-converge} can be obtained. 

Let us estimate the first of the six terms on the right-hand side of  \eqref{estimate-operator-A-2-theta-when-series-under-c-k-diverges}:
\begin{align*}
\left\|\sum_{k=1}^{\infty}   \frac{c_k}{\gamma_k}\gamma^2_k (A^2+\gamma^2_kI)^{-1} A^{2+\theta}\varphi_0 \right\|_{\mathcal{L}_{2,\gamma}}&\leq \sup_{k}\|\gamma^2_k(A^2+\gamma^2_kI)^{-1}\| \cdot \left\|\sum \limits_{k=1}^{\infty}\frac{c_k}{\gamma_k} A^{2+\theta}\varphi_0 \right\|_{\mathcal{L}_{2,\gamma}}\\
&<\frac{1}{\sqrt{2 \gamma}} \left(\sum_{k=1}^{\infty} \frac{c_k}{\gamma_k}\right) \left\| A^{2+\theta}\varphi_0 \right\|_H. 
\end{align*}
The second term on the right-hand side of inequality \eqref{estimate-operator-A-2-theta-when-series-under-c-k-diverges} is estimated as follows:
 \begin{align*}
\left\|\sum \limits_{k=1}^{\infty}   \frac{c_k}{\gamma_k} \gamma^2_k B^{-1} \mathrm{e}^{-\gamma_kt}   A^{2+\theta}\varphi_0 \right\|_{\mathcal{L}_{2,\gamma}}&\leq \sup_{k}\|\gamma^2_k(A^2+\gamma^2_kI)^{-1}\| \cdot\left\|\sum \limits_{k=1}^{\infty}   \frac{c_k}{\gamma_k}  \mathrm{e}^{-\gamma_kt} A^{2+\theta}\varphi_0 \right\|_{\mathcal{L}_{2,\gamma}}\\
&< \left(\sum \limits_{k=1}^{\infty}   \frac{c_k}{\gamma_k}\right)\left\| \mathrm{e}^{-\gamma_1t} A^{2+\theta}\varphi_0 \right\|_{\mathcal{L}_{2,\gamma}}
\end{align*} where $B=A^2+\gamma^2_kI$. Consequently, 
 \begin{align*}
\left\|\sum \limits_{k=1}^{\infty}   \frac{c_k}{\gamma_k} \gamma^2_k \left(A^2+\gamma^2_kI\right)^{-1} \mathrm{e}^{-\gamma_kt}   A^{2+\theta}\varphi_0 \right\|_{\mathcal{L}_{2,\gamma}}< \frac{1}{\sqrt{2 \gamma+\gamma_1}} \left(\sum_{k=1}^{\infty} \frac{c_k}{\gamma_k}\right) \left\| A^{2+\theta}\varphi_0 \right\|_H. 
\end{align*}

The remaining terms on the right-hand side of the inequality \eqref{estimate-operator-A-2-theta-when-series-under-c-k-diverges}  are estimated in a similar way. Indeed
 \begin{align*}
\left\|\sum \limits_{k=1}^{\infty}   \frac{c_k}{\gamma_k}\gamma_k A (A^2+\gamma^2_kI)^{-1}  A^{2+\theta}\varphi_0 \right\|_{\mathcal{L}_{2,\gamma}}&\leq \sup_{k}\|\gamma_kA(A^2+\gamma^2_kI)^{-1} \| \cdot\left\|\sum \limits_{k=1}^{\infty}   \frac{c_k}{\gamma_k}  A^{2+\theta}\varphi_0 \right\|_{\mathcal{L}_{2,\gamma}}\\
&<\frac{1}{2\sqrt{2\gamma}}\left(\sum \limits_{k=1}^{\infty}   \frac{c_k}{\gamma_k}\right)\left\| A^{2+\theta}\varphi_0 \right\|_H.
\end{align*}
\begin{align*}
\left\|\sum \limits_{k=1}^{\infty}  \frac{c_k}{\gamma_k} \gamma_k A(A^2+\gamma^2_kI)^{-1} \mathrm{e}^{-\gamma_kt}   A^{1+\theta}\varphi_1 \right\|_{\mathcal{L}_{2,\gamma}}&\leq \sup_{k}\| \gamma_k    A(A^2+\gamma^2_kI)^{-1} \| \cdot \left\|\sum \limits_{k=1}^{\infty}  \frac{c_k}{\gamma_k}  \mathrm{e}^{-\gamma_1t} A^{1+\theta}\varphi_1 \right\|_{\mathcal{L}_{2,\gamma}}\\
&\leq \frac{1}{2\sqrt{2\gamma+\gamma_1}}\left(\sum \limits_{k=1}^{\infty}  \frac{c_k}{\gamma_k}\right)\left\|  A^{1+\theta}\varphi_1 \right\|_H. 
\end{align*}
\begin{align*}
\left\|\sum \limits_{k=1}^{\infty}  \frac{c_k}{\gamma_k} \gamma_k A(A^2+\gamma^2_kI)^{-1}  A^{1+\theta}\varphi_1 \right\|_{\mathcal{L}_{2,\gamma}}&\leq \sup_{k}\|\gamma_k A(A^2+\gamma^2_kI)^{-1}\| \cdot \left\|\sum \limits_{k=1}^{\infty}  \frac{c_k}{\gamma_k} A^{1+\theta}\varphi_1 \right\|_{\mathcal{L}_{2,\gamma}}\\
&\leq \frac{1}{2\sqrt{2\gamma}}\left(\sum \limits_{k=1}^{\infty}  \frac{c_k}{\gamma_k}\right) \left\| A^{1+\theta}\varphi_1 \right\|_H. 
\end{align*}
\begin{align*}
\left\|\sum \limits_{k=1}^{\infty} \frac{c_k}{\gamma_k}  \gamma^2_k (A^2+\gamma^2_kI)^{-1} A^{1+\theta}\varphi_1 \right\|_{\mathcal{L}_{2,\gamma}}&\leq \sup_{k}\| \gamma^2_k(A^2+\gamma^2_kI)^{-1}\| \cdot \left\|\sum \limits_{k=1}^{\infty} \frac{c_k}{\gamma_k} A^{1+\theta}\varphi_1 \right\|_{\mathcal{L}_{2,\gamma}}\\
&<\frac{1}{\sqrt{2 \gamma}} \left(\sum \limits_{k=1}^{\infty} \frac{c_k}{\gamma_k}\right)\left\| A^{1+\theta}\varphi_1 \right\|.
\end{align*}  Thus, combining the estimates of six terms of the inequalities \eqref{estimate-operator-A-2-theta-when-series-under-c-k-diverges} we obtain the desired inequality \eqref{estimate-when-series-under-c-k-diverges-and-series-under-c-k-gamma-converge}. 

Now, from lemma \ref{lemma-when-series-under-c-k-diverges} and arguments similar to the point $1)$ of theorem \ref{correct-solvability-theorem-abound-interval} the following estimate is obtained 
\begin{align}\label{omega-1-when-series-under-c-k-diverges}
\|\omega_1\|_{W^2_{2} ((0, T), H)}\leq d_1(T)\left(\|A^{2+\theta}\varphi_0\|_{H}+\|A^{1+\theta}\varphi_1\|_{H}\right).
\end{align} Consequently, on the basis of estimates \eqref{estimate-under-v-with-initial-conditions}, \eqref{omega-2-when-the-series-under-c-k-diverges} and  \eqref{omega-1-when-series-under-c-k-diverges} the required estimate is established
\begin{align*}
\|u\|_{W^2_{2} ((0, T), A^2)}\leq d(T)\left(\|A^{2-\theta}f\|_{L_{2} ((0, T), H)}+\|A^{2+\theta}\varphi_0\|_{H}+\|A^{1+\theta}\varphi_1\|_{H}\right),
\end{align*} where the positive constant $d(T)$ is independent  of the vector function $f$ and vectors $\varphi_0$, $\varphi_1$.

\subsection{Proof of theorems on the location of spectra of operator-valued functions  on the left half-plane}\label{subsection-about-distribution-of zeroes-meromorph-function-ell-n}

{\bf Proof of Theorem \ref{spectr-operator-function-lying-on-the-left-right}}.  The function $\varphi(\lambda)=\lambda^2+a^2_n$, where $\lambda=x+iy$, maps the upper right quadrant $$\phi_{\pi/2}=:\{\lambda \in \mathbb{C}: 0<\arg \lambda<\pi/2\}$$ into the upper half-plane $\im \lambda>0$. In turn, the function $$\psi(\lambda)=a^{2\theta}_n\sum_{k=1}^{\infty}\frac{c_k}{\lambda+\gamma_k}$$
 maps  the angle $\phi_{\pi/2}$ into the lower half-plane $\im \lambda<0$. Therefore, the equation $\varphi(\lambda)=\psi(\lambda)$, which is equivalent to the equation $\ell_n(\lambda)=0$, has no solutions inside the angle $\phi_{\pi/2}$.  Since the function $\ell_n(\lambda)$ has real coefficients, then their non-real zeros are complex conjugate zeros. Thus, the equation $\ell_n(\lambda)=0$ has no zeros inside the lower right quadrant $\phi_{-\pi/2}=:\{\lambda \in \mathbb{C}: -\pi/2<\arg \lambda<0\}$.
 
\vspace{2.4cm}
\setlength{\unitlength}{1cm} 
\begin{picture}(1, 1)(-1, 1)
{\linethickness{0.3mm}
\put(0.5, 0){\vector(1,0){14}}}
\put(12.5,2.5){{\makebox(0,0){$\varphi(x)=x^2 +a^2_n$}}}
\put(13, 2){\vector(1,-1){0.5}}
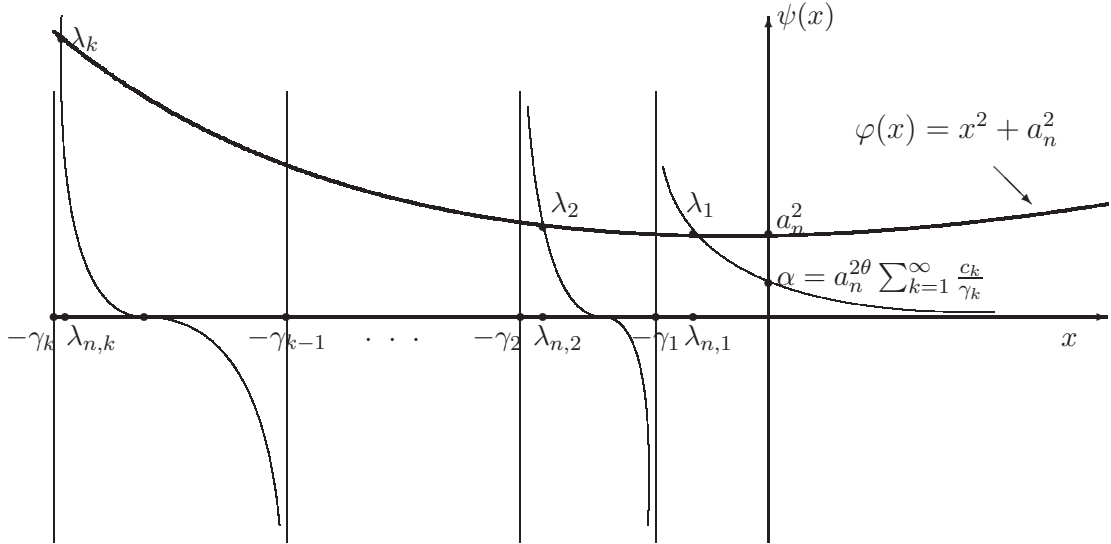
\put(7,-3.5){\makebox(0,0){{\bf Figure A}: Location of the spectrum of the operator-valued function $L(\lambda)$ on the left half-plane}}
\put(14,-0.3){\small{\makebox(0,0){$x$}}}
\put(10.3,1.3){\small{\makebox(0,0){$a^2_n$}}}
\put(11.5,0.5){\small{\makebox(0,0){$\alpha=a^{2\theta}_n\sum_{k=1}^{\infty}\frac{c_k}{\gamma_k}$}}}
{\linethickness{0.3mm}
\put(10,-3){\vector(0,1){7}}}
\put(10.5,4){\small{\makebox(0,0){$\psi(x)$}}}
{\linethickness{0.1mm}
\qbezier(1.7,0.0)(0.5,0.0)
(0.6,4)}
{\linethickness{0.1mm}
\qbezier(1.7,0.0)(3.2384,0.0)
(3.5,-2.7622)}
{\linethickness{0.34mm}
\qbezier(14.5,1.5)(5,0.0)
(0.5, 3.8)}
{\linethickness{0.1mm}
\qbezier(13,0.067)(9,0.0)
(8.6, 2)}
\multiput(8.5,0)(3.1,0){1}
{\line(0,1){3}}
\multiput(8.5,0)(3.1,0){1}
{\line(0,-1){3}}
{\linethickness{0.1mm}
\qbezier(7.8,0.0)(7.00,0.0)
(6.8, 2.8)}  
{\linethickness{0.1mm}
\qbezier(7.8,0.0)(8.5,0)
(8.4,-2.7622)}
\linethickness{.075mm}
\multiput(0.5,0)(3.1,0){3}
{\line(0,1){3}}
\put(0.2,-0.3){\small{\makebox(0,0){$-\gamma_{k}$}}}
\put(3.6,-0.3){\small{\makebox(0,0){$-\gamma_{k-1}$}}}
\put(6.4,-0.3){\small{\makebox(0,0){$-\gamma_{2}$}}}
\put(8.5,-0.3){\small{\makebox(0,0){$-\gamma_{1}$}}}
\put(5,-0.3){\makebox(0,0){$. \hspace{0.2cm}. \hspace{0.2cm}.$}}
\linethickness{.075mm}
\multiput(0.5,0)(3.1,0){3}
{\line(0,-1){3}}
\put(6.7,0.0){\circle*{0.1}}
\put(1.7,0.0){\circle*{0.1}}
\put(0.5,0.0){\circle*{0.1}}
\put(3.59,0.0){\circle*{0.1}}
\put(8.5,0.0){\circle*{0.1}}
\put(10,1.1){\circle*{0.1}}
\put(10,0.45){\circle*{0.1}}
\put(0.6,3.7){\circle*{0.1}}
\put(7,1.2){\circle*{0.1}}
\put(9,1.1){\circle*{0.1}}
\put(0.9,3.7){\small{\makebox(0,0){$\lambda_k$}}}
\put(7.3,1.5){\small{\makebox(0,0){$\lambda_2$}}}
\put(9.1,1.5){\small{\makebox(0,0){$\lambda_1$}}}
\put(0.65,0.0){\circle*{0.1}}
\put(7,0.0){\circle*{0.1}}
\put(9,0.0){\circle*{0.1}}
\put(1,-0.3){\small{\makebox(0,0){$\lambda_{n, k}$}}}
\put(7.2,-0.3){\small{\makebox(0,0){$\lambda_{n, 2}$}}}
\put(9.2,-0.3){\small{\makebox(0,0){$\lambda_{n, 1}$}}}
\end{picture}  
\vspace{5cm}

 It is also clear that when  $a^{2(1-\theta)}_n>\sum_{k=1}^{\infty}\frac{c_k}{\gamma_k}$, the equation $\varphi(x)=\psi(x)$ has no solutions, lying on the semiaxis $(0, +\infty)$ (see the figure {\bf A}).  Therefore, the equation $\varphi(\lambda)=\psi(\lambda)$ has no solutions for such $\lambda$ that $\re \lambda>0$. 

On the other hand, if  $$\sum_{k=1}^{\infty}\frac{c_k}{\gamma_k}>a^{2(1-\theta)}_n$$ is satisfied, then in the right half-plane there is at least one real zero of a meromorphic function $\ell_n(\lambda)$ (see the figure below {\bf B}). 

\vspace{2.4cm }
 \setlength{\unitlength}{1cm} 
\begin{picture}(1, 1)(-1, 1)
{\linethickness{0.3mm}
\put(0.5, 0){\vector(1,0){14}}}
\put(13,2){{\makebox(0,0){$\varphi(x)=x^2 +a^2_n$}}}
\put(13, 1.5){\vector(-1,-1){0.5}}
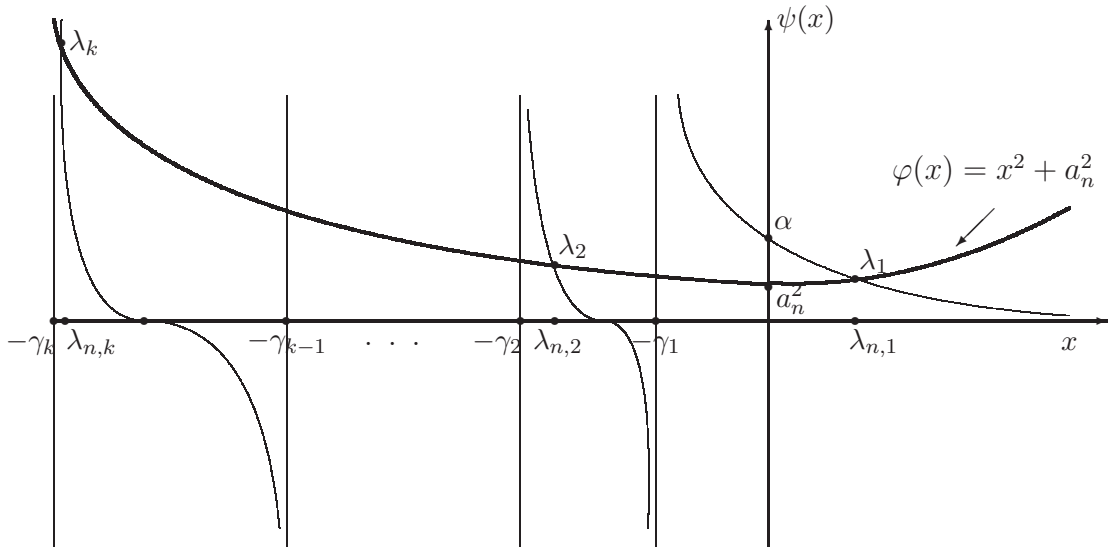
\put(6.7,-3.5){\makebox(0,0){{\bf Figure B}: Location of eigenvalues of operator-valued function $L(\lambda)$  on the right half-plane}}
\put(14,-0.3){\small{\makebox(0,0){$x$}}}
\put(10.2,1.3){\small{\makebox(0,0){$\alpha$}}}
\put(10.3,0.3){\small{\makebox(0,0){$a^2_n$}}}
{\linethickness{0.3mm}
\put(10,-3){\vector(0,1){7}}}
\put(10.5,4){\small{\makebox(0,0){$\psi(x)$}}}
{\linethickness{0.1mm}
\qbezier(1.7,0.0)(0.5,0.0)
(0.6,4)}
{\linethickness{0.1mm}
\qbezier(1.7,0.0)(3.2384,0.0)
(3.5,-2.7622)}
{\linethickness{0.34mm}
\qbezier(10,0.48)(1,1)
(0.5, 4)}
{\linethickness{0.34mm}
\qbezier(10,0.5)(12,0.4)
(14, 1.5)}
{\linethickness{0.1mm}
\qbezier(14,0.067)(9,0.4)
(8.8, 3)}
\multiput(8.5,0)(3.1,0){1}
{\line(0,1){3}}
\multiput(8.5,0)(3.1,0){1}
{\line(0,-1){3}}
{\linethickness{0.1mm}
\qbezier(7.8,0.0)(7.00,0.0)
(6.8, 2.8)}  
{\linethickness{0.1mm}
\qbezier(7.8,0.0)(8.5,0)
(8.4,-2.7622)}
\linethickness{.075mm}
\multiput(0.5,0)(3.1,0){3}
{\line(0,1){3}}
\put(0.2,-0.3){\small{\makebox(0,0){$-\gamma_{k}$}}}
\put(3.6,-0.3){\small{\makebox(0,0){$-\gamma_{k-1}$}}}
\put(6.4,-0.3){\small{\makebox(0,0){$-\gamma_{2}$}}}
\put(8.5,-0.3){\small{\makebox(0,0){$-\gamma_{1}$}}}
\put(5,-0.3){\makebox(0,0){$. \hspace{0.2cm}. \hspace{0.2cm}.$}}
\linethickness{.075mm}
\multiput(0.5,0)(3.1,0){3}
{\line(0,-1){3}}
\put(6.7,0.0){\circle*{0.1}}
\put(1.7,0.0){\circle*{0.1}}
\put(0.5,0.0){\circle*{0.1}}
\put(3.59,0.0){\circle*{0.1}}
\put(8.5,0.0){\circle*{0.1}}
\put(10,1.1){\circle*{0.1}}
\put(10,0.45){\circle*{0.1}}
\put(0.6,3.7){\circle*{0.1}}
\put(7.16,0.74){\circle*{0.1}}
\put(11.15,0.55){\circle*{0.1}}
\put(0.9,3.7){\small{\makebox(0,0){$\lambda_k$}}}
\put(7.4,1){\small{\makebox(0,0){$\lambda_2$}}}
\put(11.4,0.8){\small{\makebox(0,0){$\lambda_1$}}}
\put(0.65,0.0){\circle*{0.1}}
\put(7.16,0.0){\circle*{0.1}}
\put(11.15,0.0){\circle*{0.1}}
\put(1,-0.3){\small{\makebox(0,0){$\lambda_{n, k}$}}}
\put(7.2,-0.3){\small{\makebox(0,0){$\lambda_{n, 2}$}}}
\put(11.4,-0.3){\small{\makebox(0,0){$\lambda_{n, 1}$}}}
\end{picture}
\vspace{5cm}

In the case when 
\begin{align*}
\sum_{k=1}^{\infty}\frac{c_k}{\gamma_k}>a^{2(1-\theta)}_n, \hspace{0.2cm}n= 1, 2, \dots N_0,  \hspace{0.2cm}  a_{N_0+1}>1, 
\end{align*} on the right half-plane there is $N_0$ real zeros of a meromorphic function $\ell_n(\lambda)$. In fact, from equality $$x^2+a^{2}_j=a^{2\theta}_j\sum_{k=1}^{\infty}\frac{c_k}{x+\gamma_k},$$ for $x=0$, it is follows that
\begin{align*}
a^{2}_j=a^{2\theta}_j\sum_{k=1}^{\infty}\frac{c_k}{\gamma_k}, \hspace{0.3cm}\theta \in [0, 1] 
\end{align*} From that we have $\sum_{k=1}^{\infty}\frac{c_k}{\gamma_k}>a^{2(1-\theta)}_n$. Consequently, in the right half-plane there is $N_0$  eigenvalues of the operator-valued function $$L(\lambda)=\lambda^2 I+A^2-\widehat{K}(\lambda)A^{2\theta}$$ and the initial value problem \eqref{system 1.1}--\eqref{system 1.2} will be unstable.

Note that on the imaginary axis and at the origin, the coordinate $(0, 0)$ of complex plane $\mathbb{C}$ there are no points of the spectrum $\sigma(L)$ of operator-valued functions $L(\lambda)$. Indeed,  for $y\neq0$ and $x=0$ the following relation is valid
\begin{align*}
\im \ell_n(iy)=y\left(\sum_{k=1}^{\infty}\frac{c_k}{y^2+\gamma^2_k}\right)a^{2\theta}_n\neq 0.
\end{align*} For $x=y=0$ and on the basis of $a_1\geq 1$ and $a^{2(1-\theta)}_1>\sum_{k=1}^{\infty}\frac{c_k}{\gamma_k}$ the following equality is valid
\begin{align*}
\re \ell_n(0)=a^{2\theta}_n\left(a^{2(1-\theta)}_n-\sum_{k=1}^{\infty}\frac{c_k}{\gamma_k}\right)>0.
\end{align*}

The Theorem \ref{spectr-operator-function-lying-on-the-left-right} is proved. 


In what follows will be used the following assertions. 

We provide here a variant of the Schwarz lemma, which is known as the Schwarz--Pick lemma (see \cite[Chap. 4]{shapiro} for more details). 

{\bf  Lemma} (Schwarz--Pick \cite{shapiro} ). {\em Let $f$ be a holomorphic self-map of  the unit disc $D$. Then for every pair of points $p, q \in D$ we have}
\begin{align*}
\left|\frac{f(p)-f(q)}{1-\overline{f(p)}f(q)}\right|\leq \left|\frac{p-q}{1-\overline{p}q}\right|,
\end{align*}{\em and for all $z \in D$}, 
\begin{align*}
\frac{|f^{\prime}(z)|}{1-|f(z)|^2}\leq \frac{1}{1-|z|^2}. 
\end{align*}  There is equality for all pairs of points if and only if $f$ is a conformal automorphism of $D$.

An analogous statement on the upper half-plane $\mathbb{C}_+$ can be made as follows:

{\bf Corollary.} {\em Let $f$ be  a holomorphic self-map of  $ \mathbb{C}_+$. Then for every pair of points $p, q \in \mathbb{C}_+$} the following estimate is hold
\begin{align*}
\left|\frac{f(p)-f(q)}{\overline{f(p)}-f(q)}\right|\leq \left|\frac{p-q}{\overline{p}-q}\right|
\end{align*} and for all $z \in \mathbb{C}_+$,
\begin{align*}
 |f^{\prime}(z)|\leq \frac{\im f(z)}{\im z}.
\end{align*} 
{\bf Theorem} (Denjoy--Wolff  \cite{shapiro, shapiroBourdon}). {\em Let $f$ be an analytic function which  maps the upper half-plane $\mathbb{C}_{+}$ into itself and, suppose that $f$ is not an elliptic fractional-linear transformation. Then there exists a unique point $w \in \mathbb{C}_{+} \cup \{\infty\}$ such that the iterates $f^{*n}$ converge  to $w$ uniformly on compact subsets of $\mathbb{C}_{+}$, the angular limit $\lim\limits_{z\to w}f(z)$ exists and satisfies the equation $w=f(w)$. Moreover, the angular derivative $f^{\prime}(w)$ exists and satisfies  $f^{\prime}(w)\leq 1$}. 

\begin{remark} Angular limit means that $z$ is restricted to any angle $\epsilon<\arg (z-w)<\pi-\epsilon$, where $\epsilon>0$, if  $w 	\in \mathbb{R}$. Angular derivative is defined as $f^{\prime}(w)=\lim \limits_{z\to w}(f(z)-f(w))/(z-w)$, if  $w \in \mathbb{R}$. The transformation $f(z)$ is called an elliptic fractional-linear transformation if $f(z)$ has two fixed points in $0$ and $\infty$ (for more details, see the work \cite[chap. 2]{shapiroBourdon}). 
\end{remark}

\begin{lemma} \label{meromorfnye-funktsi-self-map-on-plus-compplex-plane} 
 Meromorphic functions $\ell_n(\lambda)$ have at most one non-real zeroes in the open upper half-plane. $\mathbb{C}_{+}$. 
\end{lemma}
 {\bf Proof of lemma \ref{meromorfnye-funktsi-self-map-on-plus-compplex-plane}}. We consider a regular branch $\varphi$ of square root, which maps the lower half-plane  $\mathbb{C}_{-}$ into the second quadrant. Then the equation
\begin{align*}
\frac{\ell_n(\lambda)}{a^2_n}=\frac{\lambda^2}{a^2_n}+1-\frac{1}{a^{2(1-\theta)}_n} \sum_{k=1}^{\infty}\frac{c_k}{\lambda+\gamma_k},
\end{align*} which is equivalent to the equation
\begin{align*}
\lambda=g(\lambda):=a_n\varphi \left(-1+\frac{1}{a^{2(1-\theta)}_n} \sum_{k=1}^{\infty}\frac{c_k}{\lambda+\gamma_k}\right),
\end{align*}i. e.,
\begin{align*}
\lambda=g(\lambda):=a_n\sqrt{-1+\frac{1}{a^{2(1-\theta)}_n} \sum_{k=1}^{\infty}\frac{c_k}{\lambda+\gamma_k}}
\end{align*}  maps the upper half-plane $\mathbb{C}_{+}$ into itself. Indeed
\begin{align*}
\im \left(\sum_{k=1}^{\infty}\frac{c_k}{\lambda+\gamma_k}\right)<0, \hspace{0.5cm}\lambda \in \mathbb{C}_{+}.
\end{align*} Therefore, by the  corollary of Schwarz--Pick lemma and  Denjoy--Wolff  theorem, the equation $\lambda=g(\lambda)$ has at most one solution in the upper half-plane $\mathbb{C}_{+}$.

\subsection{Proof of theorem on the structure of spectrum when $K(t) \in W_{1}^{1}(\mathbb{R}_+)$}\label{subsection-proof-of-theorem-about-spectral-structure-1}

{\bf Proof of the theorem \ref{theorem-about-spectral-anal-1}}. To prove the formulated theorem it will be provided 
 a series of auxiliary lemmas on the distribution of zeros of meromorphic function $\ell_n(\lambda)$ in the case when kernel $K(t)$ is written as the sum of a finite number of exponential functions.  Let us consider the function
\begin{align} \label{meromorf-function-finite}
\ell_{n,N}(\lambda):={\lambda}^2+ a^{2}_{n} \left(1- \frac{1}{a^{2(1-\theta)}_{n}}\sum_{k=1}^{N} \frac{c_k}{\lambda+\gamma_k}\right). 
\end{align}
This result is of interest, since if the conditions $a)$ and $b)$ both are satisfied, the asymptotics of zeros of meromorphic function $\ell_{n}$ can be obtained from the asymptotics of the zeros of the function $\ell_{n, N}$ by the passage to the limit when $N \to +\infty$.

\begin{lemma} \label{lemma-de-ceros-finitos} Let us consider the meromorphic function \eqref{meromorf-function-finite}. Then the zeros of the function $\ell_{n, N}$ can be written as a set of real zeros $\{\lambda_{n, k}(\theta, N)|k=1,\dots, N\}$ for which the following inequalities are satisfied
\begin{align}
   -\gamma_{k}<\lambda_{n, k}(\theta, N) <x_{n, k}(\theta, N)<-\gamma_{k-1} <\cdots <0,\label{desigualdad-ceros}
 \end{align} and such that $\lim \limits_{n\to +\infty}\lambda_{n, k}(\theta, N)=-\gamma_k$, where $x_{n, k}(\theta, N)$ are the real zeros of the function
\begin{align*}
f_{n, N}(\lambda):=1- \frac{1}{a^{2(1-\theta)}_{n}}\sum_{k=1}^{N} \frac{c_k}{\lambda+\gamma_k},  \hspace{1cm}\theta \in [0, 1]. 
\end{align*}The zeros of the function $\ell_{n, N}$ can also be written, for a sufficiently large $a_n$, as a pair of complex conjugate zeros $\lambda^{\pm}_{n}(\theta, N)$,  $\lambda^{+}_{n}(\theta, N)= \overline{\lambda^{-}_{n}}(\theta, N)$, which are asymptotically represented in the form
\begin{align} \label{asimptota-finita1}
\lambda^{\pm}_{n}(\theta, N)=-\frac{1}{2}\frac{1}{a_{n}^{2(1-\theta)}}\sum_{k=1}^{N}c_k+O\left(\frac{1}{a^{4-2\theta}_{n}}\right)\pm i\left(a_{n}+O\left(\frac{1}{a^{3-2\theta}_{n}}\right)\right), \hspace{1cm} \theta \in [0, 1]. 
\end{align} 
\end{lemma} 
\begin{lemma} \label{relations-limit-to-zeroes-complex} Suppose that the series $\sum_{k=1}^{\infty}c_k$ is convergent. Then there exists a constant $M>0$ such that for any fixed $n>M$, there is the limit
\begin{align} \label{pass-to-limit-lemma}
\lim_{N\to +\infty}\lambda^{\pm}_{n}(\theta, N)=\lambda^{\pm}_{n}(\theta). 
\end{align} 
\end{lemma}
{\bf Proof of Lemma \ref{lemma-de-ceros-finitos}.} 
It will be proved that the complex zeros of a meromorphic function $\ell_{n, N}(\lambda)$ are asymptotically represented in the form $\lambda^{\pm}_{n}(\theta, N)=\tau_n a_n \pm i a_n$, $n \in \mathbb{N}$, for $a_n \to +\infty$, where $\tau_n$ is a bounded numerical sequence. For this purpose, it is suffices to prove that the asymptotic representation $\lambda^{+}_{n}(\theta, N)$ satisfies the equation
\begin{align*}
\frac{\widehat{K}_N(\lambda^{+}_{n}(\theta, N))}{a^{2(1-\theta)}_{n}}=\frac{(\lambda^{+}_{n}(\theta, N))^2}{a^{2}_{n}}+1,
\end{align*} which is equivalent to the equation $\widehat{K}_N(\lambda^{+}_{n}(\theta, N))= a^{2(1-\theta)}_{n} \tau_n(\tau_n+2i)$. It follows that
\begin{align}\label{tau-igual-a-suma-N}
\tau_n= \frac{\widehat{K}_N(\lambda^{+}_{n}(\theta, N))}{a^{2(1-\theta)}_{n}(\tau_n+2i)}.
\end{align} We denote by 
\begin{align*}
h_n(\tau):=\frac{\widehat{K}_N( z_n)}{a^{2(1-\theta)}_{n}(\tau+2i)}, \hspace{1cm}z_n=  \tau a_n +i a_n. 
\end{align*} Then the equation \eqref{tau-igual-a-suma-N} can be rewritten in the form
$\tau_n=h_n(\tau_n)$. Hence, we get that number  $\tau_n$ is a fixed point of the mapp $\tau \to h_n(\tau)$ when $n \to +\infty$. Consequently, it suffices to prove that, when $n \to +\infty$, the mapp $\tau \to h_n(\tau)$ is a contraction mapp. Thus, the desired solution $\tau_n$ will be found as the limit of the sequence $\tau^{k}_{n}$, when $k \to +\infty$, where $\tau^{k}_{n}=h(\tau^{k-1}_{n})$, $\tau^{0}_{n}=0, \tau^{1}_{n}\not=\tau_{n}$.

\begin{remark} \label{remark-on-set-pi-minus-delta} Consider the following set $\varphi_{\pi-\delta}=\{\lambda \in \mathbb{C}: |\arg \lambda|< \pi-\delta, \delta>0\}$. Suppose that, for all $z_n \in \varphi_{\pi-\delta}$, for $|z_n|\to +\infty$, the relations $|z_n {\widehat{K}_N}'(z_n)|\to 0$ and  $|\widehat{K}_N(z_n)|\to 0$ are satisfied. Then the mapp $\tau \to h_n(\tau)$ is a contraction mapp when $a_n \to +\infty$.  Indeed, this assertion follows from the estimate
\begin{align*}
|h'_n(\tau)|=\left|\frac{\widehat{K}_N(z_n)-{\widehat{K}_N}'(z_n)(a_n\tau+2ia_n)}{a^{2(1-\theta)}_{n}(\tau+2i)^{2}}\right|\leq \frac{|{\widehat{K}_N}'(z_n)(2z_n)|}{a^{2(1-\theta)}_{n}}+\frac{|\widehat{K}_N(z_n)|}{a^{2(1-\theta)}_{n}}, 
\end{align*} for $|\tau_n|<\frac{1}{2}$. Proof of the hypothesis in the remark \ref{remark-on-set-pi-minus-delta} is given at the end of this paper. 
\end{remark}

Using the Taylor polynomial expansion in powers of $\tau_n$ we immediately obtain 
\begin{align*}
h_n(\tau_n)=-\frac{i}{2 }\frac{\widehat{K}_N(ia_n)}{a^{2(1-\theta)}_{n}}+\frac{\tau_n }{4 }\frac{\widehat{K}_N(ia_n)}{a^{2(1-\theta)}_{n}}-\frac{i \tau_n a_n}{2}\frac{\widehat{K}^{'}_N(ia_n)}{a^{2(1-\theta)}_{n}}+O(\tau^2_n).
\end{align*} It follows that
\begin{align*}
h_n(\tau_n)=-\frac{i}{2 }\frac{\widehat{K}_N(ia_n)}{a^{2(1-\theta)}_{n}}(1+O(\tau_n)), \hbox{ $n \to +\infty$}.
\end{align*} Thus, for a bounded numerical sequence $\tau_n$ the following asymptotic formula is valid: 
\begin{align}\label{formula-simpotica-de-tau-N}
\tau_n= -\frac{i}{2 }\frac{\widehat{K}_N(ia_n)}{a^{2(1-\theta)}_{n}}(1+O(\tau_n)), \mbox{ $n \to +\infty$}.
\end{align} 

Now, by the equation \eqref{formula-simpotica-de-tau-N} the following asymptotic formula will be obtained
\begin{align*}
\widehat{K}_N(ia_n)=\frac{1}{i}\left[\frac{1}{a_n}\sum_{k=1}^{N}c_k+O\left(\frac{1}{a^3_n}\right)\right]+O\left(\frac{1}{a^2_n}\right), a_n \to +\infty. 
\end{align*} Indeed, 
\begin{align*}
\widehat{K}_N(ia_n)&=\sum_{k=1}^{N}\frac{c_k}{ia_n+\gamma_k}=\sum_{k=1}^{N}\frac{c_k}{ia_n\left(1+\frac{\gamma_k}{ia_n}\right)}=\frac{1}{ia_n}\sum_{k=1}^{N}c_k \left(1+\frac{\gamma_k}{ia_n}\right)^{-1}\\
&=\frac{1}{ia_n}\sum_{k=1}^{N}c_k-\frac{1}{a^2_n}\sum_{k=1}^{N}c_k \gamma_k-\frac{1}{ia^3_n} \sum_{k=1}^{N}c_k \gamma^2_k+O\left(\frac{1}{a_n^4}\right)\\
&= \frac{1}{i}\left[\frac{1}{a_n}\sum_{k=1}^{N}c_k \right]+O\left(\frac{1}{a^2_n}\right).
\end{align*} Hence, from \eqref{formula-simpotica-de-tau-N} we obtain
\begin{align*}
\tau_n&=\left[-\frac{1}{a^{2(1-\theta)}_{n}}\left[\frac{1}{a_n}\sum_{k=1}^{N}c_k+O\left(\frac{1}{a^3_n}\right)\right]-\frac{i}{a^{2(1-\theta)}_{n}}O\left(\frac{1}{a^2_n}\right)\right](1+O(\tau_n))\\
&=-\frac{1}{2a^{3-2\theta}_{n}}\sum_{k=1}^{N}c_k+iO\left(\frac{1}{a^{4-2\theta}_n}\right). 
\end{align*} Thus, substituting $\tau_n$ in $\lambda^{\pm}_{n}(\theta, N)=\tau_n a_n \pm i a_n$, $n \in \mathbb{N}$, for  $a_n \to +\infty$, we obtain the following asymptotic representation
\begin{align} 
\lambda^{\pm}_{n}(\theta, N)=-\frac{1}{2}\frac{1}{a_{n}^{2(1-\theta)}}\sum_{k=1}^{N}c_k+O\left(\frac{1}{a^{4-2\theta}_{n}}\right)\pm i\left(a_{n}+O\left(\frac{1}{a^{3-2\theta}_{n}}\right)\right). 
\end{align} 

Graphical construction shows us that the real zeros $\{\lambda_{n, k}(\theta, N)\}_{k=1}^{N}$ of function $\ell_{n, N}(\lambda)$ satisfy the inequalities
\begin{align*}
-\gamma_{k}<\lambda_{n, k}(\theta, N)<x_{n, k}(\theta, N)<-\gamma_{k-1}, \hspace{1cm }k=1, 2, \dots, N,
\end{align*} where $x_{n, k}(\theta, N)$ are the real zeros of the function $f_{n, N}=1- \frac{1}{a^{2(1-\theta)}_{n}}\sum_{k=1}^{N} \frac{c_k}{\lambda+\gamma_k}$. It will now be proved that
\begin{align*}
 \lim\limits_{n \to +\infty}\lambda_{n, k}(\theta, N)=-\gamma_k. 
\end{align*} The real roots of the equation
\begin{align} \label{equation-with-only-real-zeroes}
\lambda^2+a^2_n-a^{2\theta}_n\sum_{k=1}^{N}\frac{c_k}{\lambda+\gamma_k}=0
\end{align} are sought in the form 
\begin{align}\label{solutions-real-zeroes}
\lambda_{n, k}(\theta, N)=-\gamma_k+\frac{c_k}{a^{2(1-\theta)}_n}+ O\left(\frac{1}{a^{2(1-\theta)}_n}\right), \hspace{0.3cm}a_n\to +\infty. 
\end{align} Then, substituting these roots in equation \eqref{equation-with-only-real-zeroes} we obtain 
\begin{align}\label{sustitution-of-solutions-in-the-original-equation}
a^{2(1-\theta)}_n+\sum_{\substack{j=1\\ j \not= k }}^{N}\frac{c_j}{\gamma_j-\gamma_k+\frac{c_k}{a^{2(1-\theta)}_n}+ O\left(\frac{1}{a^{2(1-\theta)}_n}\right)}=a^{2(1-\theta)}_n \left(1+O\left(\frac{1}{a^2_n}\right)\right).
\end{align}  Note that, for a sufficiently large $a_n$, the second term on the left-hand side of \eqref{sustitution-of-solutions-in-the-original-equation} is a bounded quantity. In fact, for any arbitrarily small $\delta>0$ there is a $N$, from which the following estimate is valid
\begin{align*}
\left|\frac{c_k}{a^{2(1-\theta)}_n}+O\left(\frac{1}{a^{2(1-\theta)}_n}\right)\right|<\delta, \hspace{0.3cm} n >N. 
\end{align*} It is not difficult to see that the following chain of inequalities holds.
{\small \begin{align*}
\left|\sum_{\substack{j=1\\ j \not= k }}^{N}\frac{c_j}{\gamma_j-\gamma_k+\frac{c_k}{a^{2(1-\theta)}_n}+ O\left(\frac{1}{a^{2(1-\theta)}_n}\right)}\right|\leq \sum_{\substack{j=1\\ j \not= k }}^{N}\frac{c_j}{\left| \gamma_j-\gamma_k-\delta \right|}\leq \sum_{\substack{j=1\\ j \not= k }}^{N}\frac{c_j}{ \gamma_j\left|1-\frac{\gamma_k+\delta}{\gamma_j} \right|}\leq \eta \sum_{\substack{j=1\\ j \not= k }}^{N}\frac{c_j}{ \gamma_j}<\eta,
\end{align*}} where $\eta=\min\left(\left|1-\frac{\gamma_k+\delta}{\gamma_{k+1}} \right|,\left|1-\frac{\gamma_k+\delta}{\gamma_k} \right| \right)$. 

We note that, for $a_n \to +\infty$, the boundedness of the second term on the left-hand side of \eqref{sustitution-of-solutions-in-the-original-equation} is obtained. Therefore, the asymptotic representation \eqref{solutions-real-zeroes} is true.  Thus, $\lim\limits_{n \to +\infty} \lambda_{n, k}(\theta, N)=-\gamma_k$,  for all $k=1, 2, \dots, N$. Lemma \ref{lemma-de-ceros-finitos} is proved.  

{\bf Proof of Lemma  \ref{relations-limit-to-zeroes-complex}}. Consider the family of equations
\begin{align}
\frac{\lambda^2 (N)}{a^2_n}+1&=\frac{1}{a^{2(1-\theta)}_n}\sum_{k=1}^{N}\frac{c_k}{\lambda (N)+\gamma_k}\label{equation-family-of-zeroes-1},\\
1&=\frac{1}{a^{2(1-\theta)}_n}\sum_{k=1}^{N}\frac{c_k}{x (N)+\gamma_k}, \label{equation-family-of-zeroes-2}
\end{align} depending on the parameter  $N \in \mathbb{N}$ for a fixed value $a_n$, where $x (N), \lambda (N) \in \mathbb{C}$ for any $N$.  In addition, we consider the equation 
\begin{align}
\frac{\lambda^2}{a^2_n}+1&=\frac{1}{a^{2(1-\theta)}_n}\sum_{k=1}^{\infty}\frac{c_k}{\lambda+\gamma_k}, \label{equation-family-of-zeroes-3}
\end{align} corresponding to the family of equations   \eqref{equation-family-of-zeroes-1}.  

Denote by $f(x):=\frac{1}{a^{2(1-\theta)}_n}\sum\limits_{k=1}^{N}\frac{c_k}{x+\gamma_k}$  the right side of equation \eqref{equation-family-of-zeroes-1} for $x\in \mathbb{R}$, for any fixed value of the parameter $N \in \mathbb{N}$. Since
\begin{align*}
\frac{df}{dx}= f^{\prime}(x)=-\frac{1}{a^{2(1-\theta)}_n}\sum_{k=1}^{N}\frac{c_k}{(x+\gamma_k)^2}<0,
\end{align*} the function $f(x)$ is a decreasing function on the set  $\mathbb{R}$ and $f(x) \to +\infty$ when $x \to -\gamma_k$, $k=1, 2, \dots, N$. Consequently, the equation \eqref{equation-family-of-zeroes-1} has $N$ real zeroes and two complex conjugate roots. The equation  \eqref{equation-family-of-zeroes-2}  has $N$ real zeroes. In turn, the equation \eqref{equation-family-of-zeroes-3} has an infinite sequence of real roots.

For any fixed value of the parameter $N \in \mathbb{N}$ let us denote by  $\lambda_{n, k}(\theta, N)$, $k=1, 2, \dots, N$, and $\lambda^{\pm}_{n}(\theta, N) = \alpha(\theta, N)\pm i\beta(\theta, N)$, $\alpha(\theta, N), \beta(\theta, N)\in \mathbb{R}$, respectively, the real and complex-conjugate roots of equation \eqref{equation-family-of-zeroes-1}. Similarly, for any fixed value of the parameter $N \in \mathbb{N}$, we denote by $x_{n, k}(\theta, N)$, $k=1, 2, \dots, N$  the real roots of the equation  \eqref{equation-family-of-zeroes-2}. In turn, by $\lambda_{n, k}(\theta)$, $k=1, 2, \dots, N$ and $\lambda^{\pm}_{n}(\theta) = \alpha_0(\theta)\pm i\beta_0(\theta)$, $\alpha_0(\theta), \beta_0(\theta)\in \mathbb{R}$ we denote, respectively, the real and complex-conjugate roots of equation  \eqref{equation-family-of-zeroes-3}. It is not difficult to see that the real roots of the equations \eqref{equation-family-of-zeroes-1} and \eqref{equation-family-of-zeroes-2} satisfy the following inequalities  
\begin{align}\label{inequality-theta-and-N-first}
-\gamma_k<\lambda_{n, k}(\theta, N)<x_{n, k}(\theta, N)<-\gamma_{k-1}, \hspace{0.5cm} k=1, 2, \dots, N. 
\end{align} In turn, the real roots of equation \eqref{equation-family-of-zeroes-3}  satisfy the following inequalities 
 \begin{align} \label{inequality-theta-and-N-second}
-\gamma_k<\lambda_{n, k}(\theta)<-\gamma_{k-1}, \hspace{0.5cm} k \in \mathbb{N}. 
\end{align} Applying  the Vieta theorem to the equation \eqref{equation-family-of-zeroes-1}  for the coefficients of powers $\lambda^{N+2}(N)$,  $\lambda^{N+1}(N)$ and free term we obtain the following relations  
 \begin{align}
 \sum_{k=1}^{N+2}\lambda_{n, k}(\theta, N)&=-\sum_{k=1}^{N} \gamma_k,  &\prod_{k=1}^{N+2}\lambda_{n, k}(\theta, N)=(-1)^{N+2} a^{2}_{n}\prod_{k=1}^{N}\gamma_k\left(1-\frac{S_N}{a^{2(1-\theta)}_{n}}\right), \label{formuli-vieta-22}
\end{align} where $S_N:=\sum \limits_{k=1}^{N} \frac{c_k}{\gamma_k}$.
In the equation corresponding to \eqref{equation-family-of-zeroes-2}, for the coefficients of powers of  $x^{N}(N)$,  $x^{N-1}(N)$ and  free term we obtain the following relations  
\begin{align}
 \sum_{k=1}^{N} x_{n, k}(\theta, N)&=\frac{\sum_{k=1}^{N} c_k }{a^{2(1-\theta)}_{n}}- \sum_{k=1}^{N} \gamma_k, &\prod_{k=1}^{N}x_{n, k}(\theta, N)=(-1)^N \prod_{k=1}^{N}\gamma_k\left(1-\frac{S_N}{a^{2(1-\theta)}_{n}}\right) \label{formuli-vieta-11} 
 \end{align} where $S_N:=\sum \limits_{k=1}^{N} \frac{c_k}{\gamma_k}$. From \eqref{formuli-vieta-22} and \eqref{formuli-vieta-11} we have
 \begin{align}
 \sum_{k=1}^{N}\lambda_{n, k}(\theta, N)+\sum_{k=1}^{N} \gamma_k&=-2\alpha (\theta, N), \label{formuli-vieta-001}\\
\prod_{k=1}^{N}\lambda_{n, k}(\theta, N)\left[\alpha^2 (\theta, N)+\beta^2 (\theta, N)\right]&=(-1)^{N+2} a^{2}_{n}\prod_{k=1}^{N}\gamma_k\left(1-\frac{S_N}{a^{2(1-\theta)}_{n}}\right), \label{formuli-vieta-002}
\end{align} 
\begin{align}
 \sum_{k=1}^{N} x_{n, k}(\theta, N)+\sum_{k=1}^{N} \gamma_k=\frac{\sum_{k=1}^{N} c_k }{a^{2(1-\theta)}_{n}},  \label{formuli-vieta-003}\\
 \prod_{k=1}^{N}x_{n, k}(\theta, N)=(-1)^N \prod_{k=1}^{N}\gamma_k\left(1-\frac{S_N}{a^{2(1-\theta)}_{n}}\right) \label{formuli-vieta-004}. 
 \end{align} From inequality  \eqref{inequality-theta-and-N-first} we obtain
 \begin{align*} 
 0<\lambda_{n, k}(\theta, N)+\gamma_k<x_{n, k}(\theta, N)+\gamma_k, \hspace{0.5cm} k=1, 2, \dots, N. 
 \end{align*} Consequently, taking into account the relations \eqref{formuli-vieta-001} and \eqref{formuli-vieta-003} we have
\begin{align}\label{sucesion-of-inequalities-when-lambda+gamma-k}
 -2\alpha (\theta, N)=  \sum_{k=1}^{N} \lambda_{n, k}(\theta, N)+\sum_{k=1}^{N} \gamma_k<\sum_{k=1}^{N} x_{n, k}(\theta, N)+\sum_{k=1}^{N} \gamma_k=\frac{\sum_{k=1}^{N} c_k }{a^{2(1-\theta)}_{n}}. 
\end{align} The sequence of sums $\sum_{k=1}^{N} c_k $ is increasing, consequently, the sequence of sums $$\sum_{k=1}^{N} x_{n, k}(\theta, N)+\sum_{k=1}^{N} \gamma_k$$ is also increasing. 
In addition, by the convergence of the series $\sum_{k=1}^{\infty} c_k$, the following estimate holds
\begin{align*}
\sum_{k=1}^{N} \left(x_{n, k}(\theta, N)+ \gamma_k\right)<\frac{\sum_{k=1}^{\infty} c_k }{a^{2(1-\theta)}_{n}}. 
\end{align*} Thus, the sequence of sums $\sum_{k=1}^{N}\left( x_{n, k}(\theta, N)+ \gamma_k\right)$ is increasing and bounded from above. Consequently,  this sequence has a limit when $N \to +\infty$. Therefore, from inequalities \eqref{sucesion-of-inequalities-when-lambda+gamma-k} we have that the sequence of sums $\sum_{k=1}^{N} \left(\lambda_{n, k}(\theta, N)+ \gamma_k\right)$ is increasing and bounded from above and, indeed, that sequence has a limit when $N\to +\infty$. 

Let us define $\psi_k(\theta, N):=\lambda_{n, k}(\theta, N)+\gamma_k$. By virtue of equality \eqref{formuli-vieta-001} there exists the limit
\begin{align} \label{relations-with-alpha-0-and-psi-k}
\lim_{N\to +\infty}\alpha(\theta, N)=-\frac{1}{2}\lim_{N\to +\infty}\sum_{k=1}^{N} \left(\lambda_{n, k}(\theta, N)+ \gamma_k\right)=-\frac{1}{2}\lim_{N\to +\infty}\sum_{k=1}^{N} \psi_k(\theta, N)=:\alpha_0(\theta). 
\end{align} 
In what follows  it will be proved that there exists the limit $\lim \limits_{N\to +\infty}\beta(\theta, N)$. In effect, from relations   \eqref{formuli-vieta-002} we obtain 
\begin{align*}
\alpha^2 (\theta, N)+\beta^2 (\theta, N)=\frac{(-1)^{N}a^{2}_{n}{\displaystyle\prod_{k=1}^{N}}\gamma_k}{{\displaystyle\prod_{k=1}^{N}}\lambda_{n, k}(\theta, N)} \left(1-\frac{S_N}{a^{2(1-\theta)}_{n}}\right). 
\end{align*} Consequently, 
\begin{align} \label{limit-on-the-right-part-of-product-formule-Vieta}
\beta^2 (\theta, N)=-\alpha^2 (\theta, N)+\frac{(-1)^N a^{2}_{n}{\displaystyle\prod_{k=1}^{N}}\gamma_k}{{\displaystyle\prod_{k=1}^{N}}\left(\psi_k(\theta, N)-\gamma_k\right)} \left(1-\frac{S_N}{a^{2(1-\theta)}_{n}}\right). 
\end{align} The right-hand side of \eqref{limit-on-the-right-part-of-product-formule-Vieta} has limit, when $N \to +\infty$. Indeed, when $N \to +\infty$ we have
\begin{align*} 
\lim_{N \to +\infty} \left(1-\frac{S_N}{a^{2(1-\theta)}_{n}}\right)= \left(1-\frac{S}{a^{2(1-\theta)}_{n}}\right),
\end{align*} where $S:=\sum\limits_{k=1}^{\infty}\frac{c_k}{\gamma_k}$. In turn, the factor 
\begin{align*} 
\frac{(-1)^N{\displaystyle\prod_{k=1}^{N}}\gamma_k}{(-1)^N{\displaystyle\prod_{k=1}^{N}}\left(\gamma_k-\psi_k(\theta, N)\right)} =\frac{(-1)^N {\displaystyle\prod_{k=1}^{N}}\gamma_k}{{\displaystyle\prod_{k=1}^{N}}\lambda_{n, k}(\theta, N)}, 
\end{align*} where $\psi_k(\theta, N)> 0$, $k\in \mathbb{N}$, can be transformed in the form
\begin{align} \label{one-divided-on-product-of-gamma_k}
\frac{1}{{\displaystyle\prod_{k=1}^{N}}\left(1-\frac{\psi_k(\theta, N)}{\gamma_k}\right)}. 
\end{align} Let us prove that the quantity \eqref{one-divided-on-product-of-gamma_k} has limit when  $N\to +\infty$. It is  not difficult to verify that the following inequality holds
\begin{align} \label{inequality-of-ln-natural-of-two-x}
-\ln (1-x)<2x, \hspace{0.5cm}x \in [0, 1/2). 
\end{align}

In turn, from the fact that $\lim\limits_{n \to +\infty}\lambda_{n, k}(\theta, N)=-\gamma_k$ (see Lemma  \ref{lemma-de-ceros-finitos}), the value $\frac{\psi_k(\theta, N)}{\gamma_k}$ tends to zero when $n \to +\infty$. Consequently, there is a constant $M>0$, such that for all $n>M$ the following inequality holds 
\begin{align*}
0<\frac{\psi_k(\theta, N)}{\gamma_k}<\frac{1}{2}. 
\end{align*} Thus, according to inequality \eqref{inequality-of-ln-natural-of-two-x} the following chain of inequalities holds
 \begin{align*} 
-\ln \prod_{k=1}^{N}\left(1-\frac{\psi_k(\theta, N)}{\gamma_k}\right)=-\sum_{k=1}^{N}\ln \left(1-\frac{\psi_k(\theta, N)}{\gamma_k}\right)\leq 2 \sum_{k=1}^{N}\frac{\psi_k(\theta, N)}{\gamma_k}\leq 2 \sum_{k=1}^{N}\psi_k(\theta, N), \hspace{0.2cm}\gamma_k\geq 1. 
\end{align*} Thus, from relation  \eqref{relations-with-alpha-0-and-psi-k} it is concluded the existence of the limit of the sequence $$-\ln \prod_{k=1}^{N}\left(1-\frac{\psi_k(\theta, N)}{\gamma_k}\right).$$ Consequently, there exists  the  limit of the sequence \eqref{one-divided-on-product-of-gamma_k}. Thus, the right-hand side of  \eqref{limit-on-the-right-part-of-product-formule-Vieta}  has limit when $N\to +\infty$. Therefore, there exists limit $\lim \limits_{N\to+\infty}\beta(N)=\beta_0(\theta)$. 
Dividing the equation \eqref{equation-family-of-zeroes-1} into its real and imaginary parts, we obtain the following relations:
\begin{align}
\frac{\alpha^2 (\theta, N)-\beta^2 (\theta, N)}{a^2_n}+1&=\frac{1}{a^{2(1-\theta)}_n}\sum_{k=1}^{N}\frac{c_k\left(\alpha (\theta, N)+\gamma_k\right)}{(\alpha (\theta, N)+\gamma_k)^2+\beta^2 (\theta, N)},\label{equation-family-of-zeroes-11}\\
\alpha (\theta, N)&=-\frac{a^{2\theta}_n}{2}\sum_{k=1}^{N}\frac{c_k}{(\alpha (\theta, N)+\gamma_k)^2+\beta^2 (\theta, N)}. \label{equation-family-of-zeroes-22}
\end{align} Passing to the limit as $N\to +\infty$ in the equation \eqref{equation-family-of-zeroes-11} and \eqref{equation-family-of-zeroes-22} we have the following equalities:
\begin{align*}
\frac{\alpha^2_0 (\theta)-\beta^2_0 (\theta)}{a^2_n}+1&=\frac{1}{a^{2(1-\theta)}_n}\sum_{k=1}^{\infty}\frac{c_k\left(\alpha_0 (\theta)+\gamma_k\right)}{(\alpha_0 (\theta)+\gamma_k)^2+\beta^2_0 (\theta)},\\
\alpha_0 (\theta)&=-\frac{a^{2\theta}_n}{2}\sum_{k=1}^{\infty}\frac{c_k}{(\alpha_0 (\theta)+\gamma_k)^2+\beta^2_0 (\theta)}.
\end{align*} Consequently, $\lambda^{\pm}_n(\theta)=\alpha_0(\theta)\pm i \beta_0(\theta)$ are the complex-conjugate roots of equation \eqref{equation-family-of-zeroes-3}. Thus, the required equality \eqref{pass-to-limit-lemma} is valid. 
Lemma \ref{relations-limit-to-zeroes-complex} is proved.

\subsection{Proof of theorem on the structure of spectrum in the case when $K(t) \notin W_{1}^{1}(\mathbb{R}_+)$, but $K(t) \in L_1(\mathbb{R}_+)$} \label{subsection-proof-of-theorem-about-spectral-structure-when-condition-not-hold}

{\bf Proof of Theorem \ref{theorem-when-condition-not-hold} }. Let us consider, as in the theorem \ref{theorem-about-spectral-anal-1}, the meromorphic function
\begin{align} \label{funtion-meromorfhic-zerousnumerables}
\frac{\ell_{n}(\lambda)}{a^{2}_{n}}= \frac{{\lambda}^2}{a^{2}_{n}}+\left(1-\frac{1}{a^{2(1-\theta)}_n}\sum_{k=1}^{\infty}\frac{c_k}{\lambda+\gamma_k}\right). 
\end{align} It will be proved that for any fixed  $n \in \mathbb{N}$  the set of all zeros of the  function $\ell_{n}(\lambda)$ is the union of a countable set of real zeros $\{\lambda_{n, k}|k \in \mathbb{N}\}$ and pairs of complex conjugate zeros  $\lambda^{\pm}_{n}$,  $\lambda^{+}_{n}= \overline{\lambda^{-}_{n}}$. Moreover, it will be shown that the sequences $\{\lambda_{n, k}|k \in \mathbb{N}\}$ satisfy the inequalities \eqref{desigualdad-ceros-finitos}. In the second part of the proof it will be established the asymptotic representations  \eqref{asimptotas-de-ceros-infinitos-1}--\eqref{asimptotas-de-ceros-infinitos-3} for complex conjugate zeros $\lambda^{\pm}_{n}$.

The main idea of the proof of the theorem \ref{theorem-when-condition-not-hold} is to replace the function $\widehat{K}(\lambda)= \sum_{k=1}^{\infty} \frac{c_k}{\lambda +\gamma_k}$ by the approximating function expressed by the integral 
\begin{align*}
h(\lambda)= \int_{1}^{\infty} \frac{\mathcal{A}}{t^{\alpha}(\lambda +\mathcal{B} t^{\beta})}dt, 
\end{align*} with subsequent estimates of the function $h(\lambda)$ in the region 
\begin{align}
\varphi_{\pi-\delta}:=\{\lambda \in
\mathbb{C}: |\arg \lambda|< \pi-\delta, \delta>0\}, 
\end{align} where $\mathcal{A}$ and $\mathcal{B}$ are some positive constants. Approximation of the function
$\widehat{K}(\lambda)$ by the function $h(\lambda)$ in the case $\theta=1$ was established in the work \cite{VR}, as well as in the monograph \cite[chap. 3]{vlasovmrau2016}.  In the case $\theta \in [0, 1]$, the proof of this assertion is given in the section \ref{proof-of-auxiliary-statements}. 

We show that the complex zeros of a meromorphic function $\ell_{n}(\lambda)$ are asymptotically represented in the form $\lambda^{\pm}_{n}=\tau_n a_n\pm i a_n $, $n \in \mathbb{N}$, for  $n \to +\infty$, where $\tau_n$ is a bounded numerical sequence. It's enough to prove that the asymptotic representation $\lambda^{+}_{n}$ satisfies the equation
\begin{align*}
\frac{\widehat{K}(\lambda^{+}_{n})}{a^{2(1-\theta)}_{n}}=\frac{(\lambda^{+}_{n})^2}{a^{2}_{n}}+1,
\end{align*} which is equivalent to equation $\widehat{K}(\lambda^{+}_{n})= a^{2(1-\theta)}_{n} \tau_n(\tau_n+2i)$. It follows that
\begin{align}\label{tau-igual-a-suma}
\tau_n= \frac{\widehat{K}(\lambda^{+}_{n})}{a^{2(1-\theta)}_{n}(\tau_n+2i)}.
\end{align} Let us denote by
\begin{align*}
h_n(\tau):=\frac{\widehat{K}(\lambda^{+}_n)}{a^{2(1-\theta)}_{n}(\tau+2i)}, \hspace{1cm}\lambda^{+}_n= \tau a_n +i a_n . 
\end{align*} Then the equation \eqref{tau-igual-a-suma} can be rewritten in the form
$\tau_n=h_n(\tau_n)$. Hence we get that number $\tau_n$ is a fixed point of the mapp $\tau \to h_n(\tau)$ for $n \to +\infty$. Consequently,  It's enough to prove that when $n \to +\infty$ the mapp $\tau \to h_n(\tau)$ is a contraction mapp. Thus, the desired solution $\tau_n$ will be found as the limit of the sequence $\tau^{k}_{n}$ when $k \to +\infty$, where $\tau^{k}_{n}=h(\tau^{k-1}_{n})$, $\tau^{0}_{n}=0, \tau^{1}_{n}\not=\tau_{n}$.

The mapp  $\tau \to h_n(\tau)$  is a contraction mapp. Indeed, this assertion follows from the estimate
\begin{align*}
|h'_n(\tau)|=\left|\frac{\widehat{K}(\lambda^{+}_{n})-{\widehat{K}}'(\lambda^{+}_{n})(a_n\tau+2ia_n)}{a^{2(1-\theta)}_{n}(\tau+2i)^{2}}\right|\leq \frac{|{\widehat{K}}'(\lambda^{+}_{n})(2\lambda^{+}_{n})|+|\widehat{K}(\lambda^{+}_{n})|}{a^{2(1-\theta)}_{n}}, 
\end{align*} where $|\tau|<\frac{1}{2}$ and  from the following lemma \ref{lemma-argument-pi-delta}, which will be proved in the section \ref{proof-of-auxiliary-statements}. 
\begin{lemma} \label{lemma-argument-pi-delta} On the set 
\begin{align*}
\varphi_{\pi-\delta}=\{\lambda \in \mathbb{C}: |\arg \lambda|< \pi-\delta, \delta>0\}, \hspace{0.3cm} \mbox{for} \hspace{0.3cm} |\lambda| \to +\infty,
\end{align*} the following relations $|\lambda {\widehat{K}}'(\lambda)|\to 0$ and  $|\widehat{K}(\lambda)|\to 0$ hold. 
 \end{lemma} 

Using the Taylor polynomial expansion  in powers $\tau_n$ we immediately obtain  
\begin{align*}
h_n(\tau_n)=-\frac{i}{2 }\frac{\hat{K}(ia_n)}{a^{2(1-\theta)}_{n}}+\frac{\tau_n}{4}\frac{\hat{K}(ia_n)}{a^{2(1-\theta)}_{n}}-\frac{i\tau_n a_n}{2}\frac{\hat{K}^{'}(ia_n)}{a^{2(1-\theta)}_{n}}+O(\tau^2_n)=-\frac{i}{2}\frac{\hat{K}(ia_n)}{a^{2(1-\theta)}_{n}}(1+O(\tau_n)), \hbox{$n \to +\infty$}.
\end{align*} 

Thus for  $\tau_n$ the following asymptotic formula holds: 
\begin{align}\label{formula-simpotica-de-tau}
\tau_n= -\frac{i}{2 }\frac{\hat{K}(ia_n)}{a^{2(1-\theta)}_{n}}(1+O(\tau_n)), \hspace{0.5cm}\mbox{$n \to +\infty$}.
\end{align} 
The following lemma is used to obtain the asymptotic formulas \eqref{asimptotas-de-ceros-infinitos-1}--\eqref{asimptotas-de-ceros-infinitos-3}. 
\begin{lemma} \label{lemma-diferen-arg-lambda} Suppose that $|\arg \lambda|< \pi-\delta, \delta>0$. Then 
\begin{align}
\widehat{K}(\lambda)&=\frac{\mathcal{A}\mathcal{B}^{r-1}}{\beta|\lambda|^{r}}\int_{0}^{\infty}\frac{dt}{t^{r}(\mathrm{e}^{i\varphi}+t)}+O\left(\frac{1}{|\lambda|}\right), &\mbox{if} \hspace{0.2cm}&\mbox{$0<r<1$},\\
\widehat{K}(\lambda)&=\frac{\mathcal{A}}{\beta}\frac{\ln\left|\frac{\lambda}{\mathcal{B}}+1\right|}{\lambda}+O\left(\frac{1}{|\lambda|}\right), &\mbox{if}  \hspace{0.2cm} &\mbox{$r=1$},
\end{align} where $\mathcal{A}$ and $\mathcal{B}$ are some positive constants. 
\end{lemma} Proof of Lemma  \ref{lemma-diferen-arg-lambda}  will be given at the end of this paper (see section \ref{proof-of-auxiliary-statements}). 

The following equality holds
\begin{align}\label{calculation-constante-D}
D=\frac{i}{2}\int_{0}^{\infty}\frac{dt}{t^{r}(i+t)}=\frac{1}{2}\frac{\pi}{\sin(\pi r)}\exp\left(-i\frac{\pi}{2}(1+r)\right).
\end{align} To verify the equality \eqref{calculation-constante-D} we need to separate the equation  \eqref{calculation-constante-D}  in its real and imaginary parts:  
\begin{align*}
\frac{i}{2}\int_{0}^{\infty}\frac{dt}{t^{r}(i+t)}=\frac{i}{2}\int_{0}^{\infty}\frac{dt}{t^{r-1}(1+t^2)} +\frac{1}{2}\int_{0}^{\infty}\frac{dt}{t^{r}(1+t^2)}. 
\end{align*} In turn, the calculation of the integrals on the right-hand side of the last equality is well known (see, for example, \cite[chap. IV]{Vollunaranovich}): 
\begin{align}
\frac{i}{2}\int_{0}^{\infty}\frac{dt}{t^{r-1}(1+t^2)}&=-\frac{i}{2}\frac{\pi \left(\res\limits_{z=i}z^{p}f(z)+\res\limits_{z=-i}z^{p}f(z)\right)}{\mathrm{e}^{\pi i p}\sin (\pi p)}=\frac{i\pi}{2}\frac{\mathrm{e}^{\pi i r}}{\sin (\pi r)}\cos \left(\frac{\pi}{2}r\right), \label{sin-exp-real-complex-1}\\
\frac{1}{2} \int_{0}^{\infty}\frac{dt}{t^{r}(1+t^2)}&= -\frac{1}{2}\frac{\pi \left(\res\limits_{z=i}z^{q}f(z)+\res\limits_{z=-i}z^{q}f(z)\right)}{\mathrm{e}^{\pi i q}\sin (\pi q)}=-\frac{\pi}{2}\frac{\mathrm{e}^{\pi i r}}{\sin (\pi r)}\sin \left(\frac{\pi}{2}r\right), \label{sin-exp-real-complex-2}
\end{align} where $p=1-r$, $f(z)=\frac{1}{1+z^2}$ and $q=-r$. Thus, from \eqref{sin-exp-real-complex-1}--\eqref{sin-exp-real-complex-2} we obtain the equality \eqref{calculation-constante-D}. 

Now, using the Lemma \ref{lemma-diferen-arg-lambda}, and the equality \eqref{calculation-constante-D}
we obtain the following asymptotic formula for the case $0<r<1$:
\begin{align*}
-\frac{i}{2}\hat{K}(ia_n)=-\frac{i}{2}\frac{\mathcal{A}\mathcal{B}^{r-1}}{\beta|a_n|^{r}}\int_{0}^{\infty}\frac{dt}{t^{r}(i+t)}+O\left(\frac{1}{a_n}\right)=-\frac{\mathcal{A}D\mathcal{B}^{1-r}}{\beta a^{r}_n} +O\left(\frac{1}{a_n}\right),
\end{align*}  Thus,
\begin{align*}
\tau_n&=  -\frac{i}{2 }\frac{\hat{K}(ia_n)}{a^{2(1-\theta)}_{n}}[1+O(\tau_n)]=\left(-\frac{\mathcal{A}D\mathcal{B}^{r-1}}{\beta a^{r+2(1-\theta)}_n} +O\left(\frac{1}{a^{1+2(1-\theta)}_n}\right)\right)[1+O(\tau_n)].
\end{align*}
Denoting $$R_1+R_2:=\left(-\frac{\mathcal{A}D_1}{\beta \mathcal{B}^{1-r}}\right)+\left(-i\frac{\mathcal{A}D_2}{\beta \mathcal{B}^{1-r}}\right),$$ where
\begin{align*}
D_1=\frac{\pi}{2}\frac{\cos\left(\frac{\pi}{2}(r+1)\right)}{\sin(\pi r)}, \hspace{1cm}D_2=-\frac{\pi}{2}\frac{\sin\left(\frac{\pi}{2}(r+1)\right)}{\sin(\pi r)},
\end{align*} we have
\begin{align*}
\tau_n &=\frac{R_1+R_2}{a^{r+2(1-\theta)}_n}+O\left(\frac{1}{a^{1+2(1-\theta)}_n}\right)+\left(\frac{R_1+R_2}{a^{r+2(1-\theta)}_n}+O\left(\frac{1}{a^{1+2(1-\theta)}_n}\right)\right) O\left(R\right)\\
&=\frac{R_1}{a^{r+2(1-\theta)}_n}+O\left(\frac{1}{a^{1+2(1-\theta)}_n}\right)+O\left(\frac{1}{a^{2r+4(1-\theta)}_n}\right)+\frac{R_2}{a^{r+2(1-\theta)}_n},
\end{align*} where $$R=\frac{R_1+R_2}{a^{r+2(1-\theta)}_n}+O\left(\frac{1}{a^{1+2(1-\theta)}_n}\right).$$ 

In the case  $0<r<1$, the following formula holds
\begin{align}
\tau_n  (\theta, r) & = \frac{R_1+R_2}{a^{r+2(1-\theta)}_n}+O\left(\frac{1}{a^{\min \{1+2(1-\theta), 2r+4(1-\theta)\}}_n}\right), \mbox{$r \in (0, \frac{1}{2})$ and $\theta \in [\frac{1}{2}, 1)$ },\label{asimpotota-tau-1}\\
\tau_n  (\theta, r)&=\frac{R_1+R_2}{a^{r+2(1-\theta)}_n}+O\left(\frac{1}{a^{1+2(1-\theta)}_n}\right), \hspace{2.1cm}\mbox{ $r \in [\frac{1}{2}, 1)$ or $\theta \in (0, \frac{1}{2})$ }.\label{asimpotota-tau-2}
\end{align} Again using the Lemma \ref{lemma-diferen-arg-lambda}, in the case $r=1$ we have
\begin{align*}
-\frac{i}{2}\hat{K}(ia_n)=-\frac{i}{2}\frac{\mathcal{A}}{\beta}\frac{\ln\left|\frac{ia_n}{\mathcal{B}}+1\right|}{ia_n}+O\left(\frac{1}{a_n}\right)=-\frac{1}{2}\frac{\mathcal{A}}{\beta}\frac{\ln a_n}{a_n}+O\left(\frac{1}{a_n}\right). 
\end{align*} Thus, we obtain the asymptotic formula 
\begin{align}\label{asimpotota-tau-3}
\tau_n (\theta, r)=-\frac{1}{2}\frac{\mathcal{A}}{\beta}\frac{\ln a_n}{a^{1+2(1-\theta)}_n}+O\left(\frac{1}{a^{1+2(1-\theta)}_n}\right). 
\end{align} Since $$\lambda^{\pm}_{n} (\theta, r)=\tau_n (\theta, r) a_n\pm ia_n,$$ then from the asymptotic formulas \eqref{asimpotota-tau-1}--\eqref{asimpotota-tau-3} we obtain the required asymptotic representations 

 \begin{align*} 
  &\lambda^{\pm}_{n}(\theta, r)=-\frac{\mathcal{B}^{r-1}\mathcal{A}D_1}{\beta a^{r+2\left(\frac{1}{2}-\theta\right)}_n}\pm i\left(a_n+\frac{\mathcal{B}^{r-1}\mathcal{A}D_2}{\beta a^{r+2\left(\frac{1}{2}-\theta\right)}_n}\right)+O\left(\frac{1}{a^{n_2(\theta, r)}_n}\right), \mbox{$r \in (0, \frac{1}{2})$ and $\theta \in [\frac{1}{2}, 1)$ }, \\
  &\lambda^{\pm}_{n}(\theta, r)=-\frac{\mathcal{B}^{r-1}\mathcal{A}D_1}{\beta a^{r+2\left(\frac{1}{2}-\theta\right)}_n}\pm i\left(a_n+\frac{\mathcal{B}^{r-1}\mathcal{A}D_2}{\beta a^{r+2\left(\frac{1}{2}-\theta\right)}_n}\right)+O\left(\frac{1}{a^{2(1-\theta)}_n}\right), \mbox{$r \in [\frac{1}{2}, 1)$ or $\theta \in (0, \frac{1}{2})$ }, \\
  &\lambda^{\pm}_{n}(\theta, r)=-\frac{1}{2}\frac{\mathcal{A}}{\beta}\frac{\ln a_n}{a^{2(1-\theta)}_n}\pm ia_n+O\left(\frac{1}{a^{2(1-\theta)}_n}\right),\hspace{5cm}\mbox{$r=1$},
\end{align*} where  $n_2(\theta, r)=\min \{2(1-\theta), 2r+3-4\theta\}$,  $r:=\frac{\alpha+\beta-1}{\beta}$, $0<\alpha \leq 1$, $\alpha +\beta>1$,  constants $\mathcal{A}>0, \mathcal{B}>0$ and constants $D_1$, $D_2$ are defined as follows
 \begin{align*} 
D_1:=\frac{\pi}{2}\frac{\cos\left(\frac{\pi}{2}(r+1)\right)}{\sin(\pi r)}, \hspace{1cm} D_2:=-\frac{\pi}{2} \frac{\sin\left(\frac{\pi}{2} (r+1)\right)}{\sin(\pi r)}. 
\end{align*} 

\subsection{Distribution of spectrum of the operator-valued function $L(\lambda)$} \label{representation-of-spectra-in-graphic-in-the-case-divergente-series}

In this section, is studied the distributions of the spectrum of the operator-valued function $L(\lambda)$  in the case, when the  kernel $K(t)$  belongs to space $L_1(\mathbb{R}_+)$, but not belongs to Sobolev space $W^{1}_{1}(\mathbb{R}_+)$. The cases  $a)$ $r=1, \theta=1$ and $ b)$ $ r= (0, 1), \theta=1$ were given in the works \cite{VR, vlasovmrau2016}. The cases when $r \in (0, 1)$ and $\theta \in [0, 1)$ are studied here.   Under the conditions of the theorem \ref{theorem-when-condition-not-hold}, from the asymptotic formulas \eqref{asimptotas-de-ceros-infinitos-1}--\eqref{asimptotas-de-ceros-infinitos-3} the  following cases are possible:

1. If $a)$ $r=1$, $\theta \in [0, 1)$  and  $b)$ $r \in (0, 1)$, $\theta \in \left(0, \frac{r+1}{2}\right)$, $\theta=1/2$,  $\im \lambda^{\pm}_{n}(\theta, r)$ asymptotically approach the imaginary axis (see figure 2), since $\re \lambda^{\pm}_{n}(\theta, r) \to -0$, for $n \to +\infty$

\vspace{1.5cm}
\setlength{\unitlength}{1cm}
\begin{picture}(1, 1)(-1, 1)     
\put(2, 0){\vector(1,0){8}}
\put(10.7,0){\makebox(0,0){$\re \lambda$}}
 \put(8,-3){\vector(0,1){6}}
\put(8.7,3){\makebox(0,0){$\im \lambda$}}
\put(7,-3.9){\makebox(0,0){{\bf Figure 2}: The structure of the spectrum in the case $K(t) \in L_1(\mathbb{R}_+)$, but $K(t) \notin W_{1}^{1}(\mathbb{R}_+)$.}}
\put(3,-1.5){\makebox(0,0){$\lambda_{n, k}(\theta)$}}
\put(3, -1.3){\vector(1,1){1}}
\put(2,0.6){\makebox(0,0){$-\infty$}}
\put(2.4, 0.3){\vector(-1,0){0.5}}
\put(2,0){\circle*{0.15}}
\put(2.1,0){\circle*{0.15}}
\put(2.2,0){\circle*{0.15}}
\put(2.3,0){\circle*{0.15}}
\put(2.4,0){\circle*{0.15}}
\put(2.8,0){\circle*{0.15}}
\put(2.9,0){\circle*{0.15}}
\put(3,0){\circle*{0.15}}
\put(3.1,0){\circle*{0.15}}
\put(3.2,0){\circle*{0.15}}
\put(3.8,0){\circle*{0.15}}
\put(3.9,0){\circle*{0.15}}
\put(4,0){\circle*{0.15}}
\put(4.1,0){\circle*{0.15}}
\put(4.2,0){\circle*{0.15}}
\put(4.8,0){\circle*{0.15}}
\put(4.9,0){\circle*{0.15}}
\put(5,0){\circle*{0.15}}
\put(5.1,0){\circle*{0.15}}
\put(5.2,0){\circle*{0.15}}
\put(4.9,2.5){\makebox(0,0){$\lambda^{\pm}_{n}(\theta, r)$}}
\put(5,2){\vector(1,-1){1}}
\put(7.3,3){\circle*{0.15}}
\put(7.1,2.2){\circle*{0.15}}
\put(6.7,1.4){\circle*{0.15}}
\put(6.2,0.8){\circle*{0.15}}
\put(5.6,0.4){\circle*{0.15}}
\put(4.9,0.15){\circle*{0.15}}
\put(4.9,-0.15){\circle*{0.15}}
\put(5.6,-0.4){\circle*{0.15}}
\put(6.2,-0.8){\circle*{0.15}}
\put(6.7,-1.4){\circle*{0.15}}
\put(7.1,-2.2){\circle*{0.15}}
\put(7.3,-3){\circle*{0.15}}
\put(8.5,1.0){\makebox(0,0){$ia_n$}}
\put(8.5,-1.0){\makebox(0,0){$-ia_n$}}
\end{picture}

\vspace{5cm}
\begin{remark} In this situation the solution of homogeneous integrodifferential equation \eqref{system 1.1} can not decay exponentially. In fact, if the solutions decayed exponentially like $\mathrm{e}^{-\alpha t}$ we would have a vertical strip $\{\lambda:\alpha<\re \lambda<0\}$ free of spectra of operator-valued function $L(\lambda)$.  
\end{remark}
2. For $r \in (0, 1)$  and $\theta \in (\frac{r+1}{2}, 1)$,  $\re \lambda^{\pm}_{n}(\theta, r) \to -\infty$ (see figure 3) when $a_n \to +\infty$, because
\begin{align}\label{figure-3-when-real-part-of-spectr-tend-to-minus-infty}
\re \lambda^{\pm}_n (\theta, r) &= -\frac{\mathcal{A}D_1}{\beta \mathcal{B}^{1-r}} a^{2\left(\theta-\frac{r+1}{2}\right)}_n+O\left(\frac{1}{a^{\min \{2(1-\theta), 2r+3-4\theta\}}_n}\right).
\end{align} 

\vspace{1.4cm}
\setlength{\unitlength}{1cm}
\begin{picture}(1, 1)(-1, 1)
\put(2, 0){\vector(1,0){8}}
\put(10.7,0){\makebox(0,0){$\re \lambda $}}
\put(7.8,-3){\vector(0,1){6}}
\put(8.7,3){\makebox(0,0){$\im \lambda$}}
\put(7,-3.7){\makebox(0,0){{\bf Figure 3}: The structure of the spectrum in the case $K(t) \in L_1(\mathbb{R}_+)$, but $K(t) \notin W_{1}^{1}(\mathbb{R}_+)$.}}
\put(1.7,3){\circle*{0.15}}
\put(2.9,1.8){\circle*{0.15}}
\put(5.5, 2.5){\makebox(0,0){$\lambda^{\pm}_n(\theta, r)$}}
\put(5.2, 2.3){\vector(-1,-1){1}}
\put(4, 1.1){\circle*{0.15}}
\put(4.9,0.7){\circle*{0.15}}
\put(5.9,0.4){\circle*{0.15}}
\put(6.8,0.2){\circle*{0.15}}
\put(6.8,-0.2){\circle*{0.15}}
\put(5.9,-0.4){\circle*{0.15}}
\put(4.9,-0.7){\circle*{0.15}}
\put(4, -1.1){\circle*{0.15}}
\put(2.9,-1.8){\circle*{0.15}}
\put(1.7,-3){\circle*{0.15}}
\put(2,0){\circle*{0.15}}
\put(2.1,0){\circle*{0.15}}
\put(2.2,0){\circle*{0.15}}
\put(2.3,0){\circle*{0.15}}
\put(2.4,0){\circle*{0.15}}
\put(2.8,0){\circle*{0.15}}
\put(2.9,0){\circle*{0.15}}
\put(3,0){\circle*{0.15}}
\put(3.1,0){\circle*{0.15}}
\put(3.2,0){\circle*{0.15}}
\put(3.8,0){\circle*{0.15}}
\put(3.9,0){\circle*{0.15}}
\put(4,0){\circle*{0.15}}
\put(4.1,0){\circle*{0.15}}
\put(4.2,0){\circle*{0.15}}
\put(4.8,0){\circle*{0.15}}
\put(4.9,0){\circle*{0.15}}
\put(5,0){\circle*{0.15}}
\put(5.1,0){\circle*{0.15}}
\put(5.2,0){\circle*{0.15}}
\put(5.8,0){\circle*{0.15}}
\put(5.9,0){\circle*{0.15}}
\put(6,0){\circle*{0.15}}
\put(6.1,0){\circle*{0.15}}
\put(6.2,0){\circle*{0.15}}
\put(6.8,0){\circle*{0.15}}
\put(6.9,0){\circle*{0.15}}
\put(7,0){\circle*{0.15}}
\put(7.1,0){\circle*{0.15}}
\put(7.2,0){\circle*{0.15}}
\put(8.5,1){\makebox(0,0){$ia_n$}}
\put(8.5,-1){\makebox(0,0){$-ia_n$}}
\put(2,0.6){\makebox(0,0){$-\infty$}}
\put(2,-1.3){\makebox(0,0){$\lambda_{n, k}(\theta)$}}
\put(2, -1.1){\vector(1,1){1}}
\put(2.4, 0.3){\vector(-1,0){0.5}}
\end{picture}

\vspace{5.5cm}

3.  For $r \in (0, 1)$ and $\theta= \frac{r+1}{2}$, the following asymptotic formula is obtained 
\begin{align} \label{asymptotic-formula-when-real-part-equal-minus-constant}
\re \lambda^{\pm}_n (r) = -\frac{\mathcal{A} D_1}{\beta \mathcal{B}^{1-r}}+O\left(\frac{1}{a^{1-r}_{n}}\right),
\end{align} which depends only on $r$. Consequently, when $a_n \to +\infty$,  $ \re \lambda^{\pm}_{n}(\theta, r) \to \vartheta$, where $\vartheta:=-\frac{\mathcal{A}D_1}{\beta \mathcal{B}^{1-r}}$. We note that when $\mathcal{A}\geq \mathcal{B}$, the spectrum is shown in the figure 4, but if $\mathcal{A}<\mathcal{B}$, then  the spectrum is shown in the figure 5. 

\vspace{2cm}

\setlength{\unitlength}{1cm}
\begin{picture}(1, 1)(-1, 1)     
\begin{centering}
\put(2, 0){\vector(1,0){8}}
\put(10.7,0){\makebox(0,0){ $\re\lambda$}}
 \put(8,-3){\vector(0,1){6}}
\put(8.7,3){\makebox(0,0){$\im \lambda$}}
\put(7,-4){\makebox(0,0){{\bf Figure 4}: The structure of the spectrum in the case $K(t) \in L_1(\mathbb{R}_+)$, but $K(t) \notin W_{1}^{1}(\mathbb{R}_+)$.}}
\put(2.8,-1.5){\makebox(0,0){$\lambda_{n, k}(\theta)$}}
\put(2.8, -1.3){\vector(1,1){1}}
\put(2,0.6){\makebox(0,0){$-\infty$}}
\put(2.4, 0.3){\vector(-1,0){0.5}}
\put(2,0){\circle*{0.15}}
\put(2.1,0){\circle*{0.15}}
\put(2.2,0){\circle*{0.15}}
\put(2.3,0){\circle*{0.15}}
\put(2.4,0){\circle*{0.15}}
\put(2.8,0){\circle*{0.15}}
\put(2.9,0){\circle*{0.15}}
\put(3,0){\circle*{0.15}}
\put(3.1,0){\circle*{0.15}}
\put(3.2,0){\circle*{0.15}}
\put(3.7,0){\circle*{0.15}}
\put(3.8,0){\circle*{0.15}}
\put(3.9,0){\circle*{0.15}}
\put(4,0){\circle*{0.15}}
\put(4.1,0){\circle*{0.15}}
\put(4.5,0){\circle*{0.15}}
\put(4.6,0){\circle*{0.15}}
\put(4.7,0){\circle*{0.15}}
\put(4.8,0){\circle*{0.15}}
\put(4.9,0){\circle*{0.15}}
\put(4.9,2.5){\makebox(0,0){$\lambda^{\pm}_{n}(\theta, r)$}}
\put(5,2){\vector(1,-1){1}}
\put(7.3,3){\circle*{0.15}}
\put(7.1,2.2){\circle*{0.15}}
\put(6.7,1.4){\circle*{0.15}}
\put(6.2,0.8){\circle*{0.15}}
\put(5.6,0.4){\circle*{0.15}}
\put(4.9,0.15){\circle*{0.15}}
\put(4.9,-0.15){\circle*{0.15}}
\put(5.6,-0.4){\circle*{0.15}}
\put(6.2,-0.8){\circle*{0.15}}
\put(6.7,-1.4){\circle*{0.15}}
\put(7.1,-2.2){\circle*{0.15}}
\put(7.3,-3){\circle*{0.15}}
\put(8.5,2){\makebox(0,0){$ia_n$}}
\put(8.5,-2){\makebox(0,0){$-ia_n$}}
\put(8.7,1){\makebox(0,0){$\vartheta$}} 
\put(8.5, 1){\vector(-1,-1){1}}
\linethickness{0.15mm}
\multiput(7.5,0)(3.1,0){1}{\line(0,1){3}}
\linethickness{0.15mm}
\multiput(7.5,0)(3.1,0){1}{\line(0,-1){3}}
\end{centering}
\end{picture}

\vspace{5.5cm}
It is pertinent to mention that the structure of the complex spectrum of operator-valued functions $L(\lambda)$, considered here, differs significantly from the structure of the complex spectrum of the operator-valued function $L(\lambda)$ of work \cite{VR}. Indeed, from the asymptotic formulas \eqref{asimptotas-de-ceros-infinitos-1}--\eqref{asimptotas-de-ceros-infinitos-3} for $\theta=1$ we have $\re \lambda^{\pm}_{n}(\theta, r) \to -\infty$,  when $n \to +\infty$ (for more details see work\cite{VR}). In the case $\theta \in [0, 1)$ we have  $\re \lambda^{\pm}_{n}(\theta, r) \to -\infty$,  when $n \to +\infty$,  (see figure 3) or $ \re \lambda^{\pm}_{n}(\theta, r) \to 0$  (see figure 2) or  $\re \lambda^{\pm}_{n}(\theta, r)$ tend to a negative constant (see figures 4 and 5).  

\vspace{2cm}

\setlength{\unitlength}{1cm}
\begin{picture}(1, 1)(-1, 1)
\put(2, 0){\vector(1,0){8}}
\put(10.7,0){\makebox(0,0){$\re \lambda $}}
\put(7,-3.9){\makebox(0,0){{\bf Figure 5}: The structure of the spectrum in the case $K(t) \in L_1(\mathbb{R}_+)$, but $K(t) \notin W_{1}^{1}(\mathbb{R}_+)$.}}
\put(7.8,-3){\vector(0,1){6}}
\put(8.7,3){\makebox(0,0){$\im \lambda$}}
\put(3.1,2.8){\circle*{0.15}}
\put(3.3,2.3){\circle*{0.15}}
\put(3.5,1.9){\circle*{0.15}}
\put(4.4,2){\makebox(0,0){$\lambda^{\pm}_n(\theta, r)$}}
\put(4, 1.3){\circle*{0.15}}
\put(4.9,0.6){\circle*{0.15}}
\put(5.9,0.25){\circle*{0.15}}
\put(6.8,0.15){\circle*{0.15}}
\put(6.8,-0.15){\circle*{0.15}}
\put(5.9,-0.25){\circle*{0.15}}
\put(4.9,-0.6){\circle*{0.15}}
\put(4, -1.3){\circle*{0.15}}
\put(3.5,-1.9){\circle*{0.15}}
\put(3.3,-2.3){\circle*{0.15}}
\put(3.1,-2.8){\circle*{0.15}}
\put(4,-0.4){\makebox(0,0){$\lambda_{n, k}(\theta)$}}
\put(3.8,0){\circle*{0.15}}
\put(3.9,0){\circle*{0.15}}
\put(4,0){\circle*{0.15}}
\put(4.1,0){\circle*{0.15}}
\put(4.2,0){\circle*{0.15}}
\put(4.8,0){\circle*{0.15}}
\put(4.9,0){\circle*{0.15}}
\put(5,0){\circle*{0.15}}
\put(5.1,0){\circle*{0.15}}
\put(5.2,0){\circle*{0.15}}
\put(5.8,0){\circle*{0.15}}
\put(5.9,0){\circle*{0.15}}
\put(6,0){\circle*{0.15}}
\put(6.1,0){\circle*{0.15}}
\put(6.2,0){\circle*{0.15}}
\put(6.9,0){\circle*{0.15}}
\put(7,0){\circle*{0.15}}
\put(7.1,0){\circle*{0.15}}
\put(7.2,0){\circle*{0.15}}
\put(7.3,0){\circle*{0.15}}
\put(8.5,1.0){\makebox(0,0){$ia_n$}}
\put(8.5,-1.0){\makebox(0,0){$-ia_n$}}
\linethickness{0.15mm}
\multiput(2.8,0)(3.1,0){1}{\line(0,1){3}}
\linethickness{0.15mm}
\multiput(2.8,0)(3.1,0){1}{\line(0,-1){3}}
\put(1.8,1.5){\makebox(0,0){$\vartheta$}} 
\put(1.75, 1.1){\vector(1,-1){1}}
\end{picture}
 
\vspace{5.5cm}

The presence of the parameter $\theta \in [0, 1)$ complicates the structure of the non-real spectrum of the operator-valued function $L(\lambda)$. For instance, in the cases $1)$ and $3)$ (see figures 2, 4 and 5) the structure of the non-real spectrum of the operator-valued function $L(\lambda)$ is close to the spectrum of the wave equation, while in the case $2)$ (see figure 3) the structure of the non-real spectrum of $L(\lambda)$ is close to the spectrum of an abstract parabolic equation.

\section{Proof of auxiliary statements}\label{proof-of-auxiliary-statements}
\begin{proposition} \label{propsition-of difference-arg-lambda-and-real-t} Suppose that $|\arg \lambda|<\pi-\delta$, with $\delta>0$. Then for  $\lambda\in \mathbb{C}$ and $t\in \mathbb{R}$ we have
\begin{align*}
|\lambda+t|^2 \approx |\lambda|^2 + |t|^2. 
\end{align*}
\end{proposition}
Let us denote by $\varphi:=\arg \lambda$. Then on the one hand we have
\begin{align*} 
|\lambda+t|^2=|\lambda|^2+ |t|^2+2|\lambda||t|\cos \varphi\leq 2\left(|\lambda|^2+ |t|^2\right). 
\end{align*} On the other hand, the following relations are obtained
\begin{align*} 
|\lambda+t|^2=|\lambda|^2+ |t|^2+2|\lambda||t|\cos \varphi \geq|\lambda|^2+ |t|^2-2|\lambda||t|\cos \delta\geq(1-\cos \delta)(|\lambda|^2+ |t|^2). 
\end{align*}
\begin{lemma} \label{lemma-about-difference-of-kernel-and-approximation-function-h-from-lambda} Let $h(\lambda)= \int_{1}^{\infty} \frac{\mathcal{A}}{x^{\alpha}(\lambda +\mathcal{B} x^{\beta})}dx$ and suppose that $|\arg \lambda|<\pi-\delta$, where $\delta>0$. Then the following estimate holds
\begin{align*}
\left|\widehat{K}(\lambda)-h(\lambda)\right|\leq  \frac{\operatorname{const}}{|\lambda|}.
\end{align*}
\end{lemma}
{\bf Proof.} The sequences $c_k$ and $\gamma_k$ are represented in the form 
\begin{align*}
c_k=\frac{\mathcal{A}}{k^{\alpha}}+\varphi(k), \hspace{0.5cm} \gamma_k=\mathcal{B}k^{\beta}+\phi(k), 
\end{align*} where  $\mathcal{A}>0, \mathcal{B}>0$, $0<\alpha \leq 1$, $\alpha +\beta>1$ and $\varphi(k)=O\left(\frac{1}{k^{\alpha+1}}\right)$, $\phi(k)=O\left(k^{\beta-1}\right)$, for $k\to +\infty$. Now, let us estimate the modulus of the difference
 \begin{align}
\left|\widehat{K}(\lambda)-h(\lambda)\right|&= \left|\sum_{k=1}^{\infty}\int_{k}^{k+1}\left(\frac{\mathcal{A}}{x^{\alpha}(\lambda +\mathcal{B} x^{\beta})}-\frac{c_k}{\lambda+\gamma_k}\right)dx\right| \leq \sum_{k=1}^{\infty}\int_{k}^{k+1}\left|\frac{\mathcal{A}}{x^{\alpha}(\lambda +\mathcal{B} x^{\beta})}-\frac{c_k}{\lambda+\gamma_k}\right|dx\nonumber\\
&\leq \sum_{k=1}^{\infty}\int\limits _{k}^{k+1}\left[\left|\frac{\frac{\mathcal{A}}{x^{\alpha}}}{\lambda +\mathcal{B} x^{\beta}}-\frac{\frac{\mathcal{A}}{k^{\alpha}}}{\lambda +\mathcal{B} k^{\beta}}\right|+\left|\frac{\frac{\mathcal{A}}{k^{\alpha}}+\varphi(k)}{\lambda+\mathcal{B}k^{\beta}+\phi(k)}-\frac{\frac{\mathcal{A}}{k^{\alpha}}}{\lambda +\mathcal{B} k^{\beta}}\right|\right] dx. \label{difference-estimate-of-kernel-and-approximation-function-h-from-lambda}
\end{align} We denote by $r(x):=\frac{1}{x^{\alpha}(\lambda+\mathcal{B}x^{\beta})}$. We estimate the first term in formula \eqref{difference-estimate-of-kernel-and-approximation-function-h-from-lambda}. It is observed, at first, that for all $x \in [k, k+1]$ the following estimate is true 
\begin{align*}
|r(x)-r(k)|\leq \max_{x\in  [k, k+1] }|r^{\prime}(x)|. 
\end{align*} According to proposition \ref{propsition-of difference-arg-lambda-and-real-t}, the following estimate is obtained
\begin{align*}
|r^{\prime}(x)|&= \left|\frac{\alpha x^{\alpha-1}(\lambda+\mathcal{B}x^{\beta})+\beta\mathcal{B}x^{\alpha+\beta-1}}{x^{2\alpha}(\lambda+\mathcal{B}x^{\beta})^2}\right|=\left|\frac{\alpha (k+1)^{\alpha-1}\lambda+(\alpha+\beta) \mathcal{B}(k+1)^{\alpha+\beta-1}}{k^{2\alpha}(\lambda+\mathcal{B}k^{\beta})^2}\right|\\
&\leq \frac{(k+1)^{\alpha-1}\alpha|\lambda|+(\alpha+\beta) \mathcal{B}(k+1)^{\alpha+\beta-1}}{k^{2\alpha}\left|\lambda+\mathcal{B}k^{\beta}\right|^2}\\
&\lesssim \frac{1}{k^{\alpha+\beta+1}} \left(1+\frac{1}{k}\right)^{\alpha}\frac{1}{|\lambda|}\left(\frac{\frac{\alpha}{k^{\beta}}+\frac{(\alpha+\beta) \mathcal{B}\left( 1+\frac{1}{k}\right)^{\beta}}{|\lambda|}}{\frac{1}{k^{2\beta}}+\frac{\mathcal{B}}{|\lambda|^2}}\right)\lesssim \frac{1}{k^{\alpha+\beta+1}} \frac{1}{|\lambda|}.
\end{align*} Thus, the first modulus in the formula \eqref{difference-estimate-of-kernel-and-approximation-function-h-from-lambda} is established
\begin{align}\label{estimate-firts-term-by-sum-k-potential-alpha-plus-one}
\sum_{k=1}^{\infty}\int\limits _{k}^{k+1}\left|\frac{\frac{\mathcal{A}}{x^{\alpha}}}{\lambda +\mathcal{B} x^{\beta}}-\frac{\frac{\mathcal{A}}{k^{\alpha}}}{\lambda +\mathcal{B} k^{\beta}}\right|dx\lesssim  \frac{1}{|\lambda|}\sum_{k=1}^{\infty} \frac{1}{k^{\alpha+\beta+1}}. 
\end{align} The second modulus in the formula \eqref{difference-estimate-of-kernel-and-approximation-function-h-from-lambda} is established by the proposition \ref{propsition-of difference-arg-lambda-and-real-t}: 
{\small \begin{align*}
\left|\frac{\frac{\mathcal{A}}{k^{\alpha}}+\varphi(k)}{\lambda+\mathcal{B}k^{\beta}+\phi(k)}-\frac{\frac{\mathcal{A}}{k^{\alpha}}}{\lambda +\mathcal{B} k^{\beta}}\right|^2&= \left|\frac{\varphi(k) \lambda+\mathcal{B}\varphi(k) k^{\beta}-\frac{\mathcal{A}}{k^{\alpha}}\phi(k)}{(\lambda+\mathcal{B}k^{\beta}+\phi(k) )(\lambda +\mathcal{B} k^{\beta})}\right|^2
&\lesssim \frac{\frac{|\lambda|^2}{k^{2(\alpha+1)}}\left[O\left(\frac{1}{k^{4\beta}}\right) +\frac{\left(\mathcal{B}-\mathcal{A}\right)^2}{|\lambda|^2} O\left(\frac{1}{k^{2\beta}}\right)\right]}{\left(\frac{\left|\lambda\right|^2}{k^{2\beta}}+\left|\mathcal{B}+O\left(\frac{1}{k}\right) \right|^2\right)\left(\frac{\left|\lambda\right|^2}{k^{2\beta}} +\mathcal{B}^2 \right)}.
\end{align*}} Hence it follows 
\begin{align*}
\left|\frac{\frac{\mathcal{A}}{k^{\alpha}}+\varphi(k)}{\lambda+\mathcal{B}k^{\beta}+\phi(k)}-\frac{\frac{\mathcal{A}}{k^{\alpha}}}{\lambda +\mathcal{B} k^{\beta}}\right|^2&\lesssim \frac{1}{k^{2(\alpha+1)}}\frac{1}{|\lambda|^2}\frac{\left[O\left(\frac{1}{k^{4\beta}}\right) +\frac{\left(\mathcal{B}-\mathcal{A}\right)^2}{|\lambda|^2} O\left(\frac{1}{k^{2\beta}}\right)\right]}{\left(\frac{1}{k^{2\beta}}+\frac{\left|\mathcal{B}+O\left(\frac{1}{k}\right) \right|^2}{|\lambda|^2}\right)\left(\frac{1}{k^{2\beta}} +\frac{\mathcal{B}^2 }{\left|\lambda\right|^2}\right)}
\lesssim \frac{1}{k^{2(\alpha+1)}}\frac{1}{|\lambda|^2}. 
\end{align*} Consequently,
\begin{align} \label{estimate-second-term-by-sum-k-potential-alpha-plus-one}
\sum_{k=1}^{\infty}\int\limits _{k}^{k+1}\left|\frac{\frac{\mathcal{A}}{k^{\alpha}}+\varphi(k)}{\lambda+\mathcal{B}k^{\beta}+\phi(k)}-\frac{\frac{\mathcal{A}}{k^{\alpha}}}{\lambda +\mathcal{B} k^{\beta}}\right| dx\lesssim\frac{1}{|\lambda|}\sum_{k=1}^{\infty} \frac{1}{k^{\alpha+1}}. 
\end{align} Therefore, from estimates \eqref{estimate-firts-term-by-sum-k-potential-alpha-plus-one} and \eqref{estimate-second-term-by-sum-k-potential-alpha-plus-one} the following inequalities are established
\begin{align}
\left|\widehat{K}(\lambda)-h(\lambda)\right|\lesssim \frac{1}{|\lambda|}\sum_{k=1}^{\infty} \frac{1}{k^{\alpha+\beta+1}}+\frac{1}{|\lambda|}\sum_{k=1}^{\infty} \frac{1}{k^{\alpha+1}}\leq  \frac{\operatorname{const}}{|\lambda|},
\end{align}  where  $0<\alpha \leq 1$, $\alpha +\beta>1$. Lemma \ref{lemma-about-difference-of-kernel-and-approximation-function-h-from-lambda} is proved. 

{\bf Proof of Lemma \ref{lemma-diferen-arg-lambda}.} Let $0<r<1$. Making a change of variables $\mathcal{B}x^{\beta}=|\lambda| t$ in the integral $h(\lambda)$ immediately is obtained
 \begin{align*}
h(\lambda)=\int\limits_{1}^{\infty}\frac{\mathcal{A}dx}{x^{\alpha}(\lambda+\mathcal{B}x^{\beta})}= \frac{1}{\beta |\lambda|^r}\int\limits_{\frac{\mathcal{B}}{|\lambda|}}^{\infty}\frac{\mathcal{A}\mathcal{B}^{r-1}dt}{t^{r}(\mathrm{e}^{i\varphi}+t)}=\frac{1}{\beta |\lambda|^r}\int\limits_{\frac{\mathcal{B}}{|\lambda|}}^{0}\frac{\mathcal{A}\mathcal{B}^{r-1}dt}{t^{r}(\mathrm{e}^{i\varphi}+t)} +\frac{1}{\beta |\lambda|^r}\int\limits_{0}^{\infty}\frac{\mathcal{A}\mathcal{B}^{r-1}dt}{t^{r}(\mathrm{e}^{i\varphi}+t)},
\end{align*} where $r=\frac{\alpha+\beta-1}{\beta}$, $\mathcal{A}>0, \mathcal{B}>0$, $0<\alpha \leq 1$, $\alpha +\beta>1$.  The following integral can be computed using basic method of integration
\begin{align*}
\int\limits_{\frac{\mathcal{B}}{|\lambda|}}^{0}&\frac{\mathcal{A}\mathcal{B}^{r-1}dt}{t^{r}(\mathrm{e}^{i\varphi}+t)}= \mathrm{e}^{-i\varphi}\int\limits_{\frac{\mathcal{B}}{|\lambda|}}^{0}\frac{\mathcal{A}\mathcal{B}^{r-1}dt}{t^{r}(1+t\mathrm{e}^{-i\varphi})}= \mathrm{e}^{-i\varphi}\int\limits_{\frac{\mathcal{B}}{|\lambda|}}^{0}\frac{\mathcal{A}\mathcal{B}^{r-1}}{t^{r}}\left[1-\frac{t}{\mathrm{e}^{i\varphi}}+\frac{t^2}{\mathrm{e}^{2i\varphi}}+\cdots\right]dt\\
&=-\frac{\mathrm{e}^{-i\varphi}\mathcal{B}^{1-r}}{(1-r)|\lambda|^{1-r}}+\frac{\mathrm{e}^{-2i\varphi}\mathcal{B}^{2-r}}{(2-r)|\lambda|^{2-r}}-\frac{\mathrm{e}^{-3i\varphi}\mathcal{B}^{3-r}}{(3-r)|\lambda|^{3-r}}+\cdots=-\frac{\mathrm{e}^{-i\varphi}\mathcal{B}^{1-r}}{(1-r)|\lambda|^{1-r}}+O\left(\frac{1}{|\lambda|^{2-r}}\right).
\end{align*} Consequently, 
\begin{align}\label{calculation-integral-of-difference-of-kernel-and-approximation-function}
h(\lambda)&=\frac{\beta^{-1}}{ |\lambda|^r}\left(-\frac{\mathrm{e}^{-i\varphi}\mathcal{B}^{1-r}}{(1-r)|\lambda|^{1-r}}+O\left(\frac{1}{|\lambda|^{2-r}}\right)+\int\limits_{0}^{\infty}\frac{\mathcal{A}\mathcal{B}^{r-1}dt}{t^{r}(\mathrm{e}^{i\varphi}+t)}\right)\nonumber\\
&=\frac{\beta^{-1}}{ |\lambda|^r}\int\limits_{0}^{\infty}\frac{\mathcal{A}\mathcal{B}^{r-1}dt}{t^{r}(\mathrm{e}^{i\varphi}+t)}-\frac{1}{\beta}\frac{\mathrm{e}^{-i\varphi}\mathcal{B}^{1-r}}{(1-r)|\lambda|}+O\left(\frac{1}{|\lambda|^{2}}\right)= \frac{\mathcal{A}\mathcal{B}^{r-1}}{\beta |\lambda|^r}\int\limits_{0}^{\infty}\frac{t^{-r}dt}{(\mathrm{e}^{i\varphi}+t)}+O\left(\frac{1}{|\lambda|}\right).
\end{align}  
In the case $r=1$, it is enough to make a change of variables $\mathcal{B}x^{\beta}=t$ in the integral $h(\lambda)$: 

\begin{align} \label{calculation-integral-of difference-of-kernel-and-approximation-function-when-r-equal-one}
h(\lambda)=\int\limits_{1}^{\infty}\frac{\mathcal{A}dx}{x^{\alpha}(\lambda+\mathcal{B}x^{\beta})}= \frac{\mathcal{A}}{\beta}\int\limits_{\mathcal{B}}^{\infty}\frac{dt}{t(\lambda+t)}= \frac{\mathcal{A}}{\beta}\frac{\ln \left(\frac{\lambda}{\mathcal{B}}+1\right)}{\lambda}. 
\end{align} From Lemma \eqref{lemma-about-difference-of-kernel-and-approximation-function-h-from-lambda} and estimates  \eqref{calculation-integral-of-difference-of-kernel-and-approximation-function} and \eqref{calculation-integral-of difference-of-kernel-and-approximation-function-when-r-equal-one} we obtain the proof of the Lemma \ref{lemma-diferen-arg-lambda}. 
\begin{corollary} \label{corollary-about-estimates-complex-and-real-part-of-function-h-from-i-y} The following relation is immediately fulfilled 
\begin{align*} 
|h(iy)|=\frac{\mathcal{A}\mathcal{B}^{r-1}}{\beta |y|^r}\left|\int_{0}^{\infty}\frac{dt}{t^r(i+t)}\right|<+\infty, \hspace{0.3cm} y\in \mathbb{R}. 
\end{align*} 
\end{corollary} Indeed,
\begin{align*} 
|\im h(iy)|=\frac{\mathcal{A}\mathcal{B}^{r-1}}{\beta |y|^r}\left|\int_{0}^{\infty}\im \frac{dt}{t^r(i+t)}\right|=\frac{\mathcal{A}\mathcal{B}^{r-1}}{\beta |y|^r}\left|\int_{0}^{\infty}\frac{dt}{t^r(t^2+1)}\right|<+\infty, \\
|\re h(iy)|=\frac{\mathcal{A}\mathcal{B}^{r-1}}{\beta |y|^r}\left|\int_{0}^{\infty} \re \frac{dt}{t^r(i+t)}\right|=\frac{\mathcal{A}\mathcal{B}^{r-1}}{\beta |y|^r}\left|\int_{0}^{\infty} \frac{dt}{t^{r-1}(t^2+1)}\right|<+\infty. 
\end{align*}

{\bf Proof of Lemma \ref{lemma-argument-pi-delta}.} From the Corollary \ref{corollary-about-estimates-complex-and-real-part-of-function-h-from-i-y}, the following inequalities are established  immediately
\begin{align*}
|\widehat{K}(\lambda)|&\leq \left|\frac{\mathcal{A}\mathcal{B}^{r-1}}{\beta|\lambda|^{r}}\int_{0}^{\infty}\frac{dt}{t^{r}(\mathrm{e}^{i\varphi}+t)}\right|+\left| O\left(\frac{1}{|\lambda|}\right)\right|\\
&\lesssim \frac{\mathcal{A}\mathcal{B}^{r-1}}{\beta|\lambda|^{r}}\left(\left|\int_{0}^{\infty}\frac{dt}{t^r(t^2+1)}\right|+\left|\int_{0}^{\infty} \frac{dt}{t^{r-1}(t^2+1)}\right|\right) +\left| O\left(\frac{1}{|\lambda|}\right)\right|.
\end{align*} Thus, $|\widehat{K}(\lambda)| \to 0$ when $\lambda \to +\infty$. Additionally, by the Proposition \ref{propsition-of difference-arg-lambda-and-real-t} the following estimate is obtained
\begin{align} \label{estimate-of-derivate-kernel-function-by-lambda}
|\lambda\widehat{K}^{'}(\lambda)|\leq \sum_{k=1}^{\infty}\frac{c_k |\lambda| }{\left|\lambda+\gamma_k\right|^2}\lesssim \sum_{k=1}^{\infty}\frac{c_k |\lambda| }{\left|\lambda\right|^2+\left|\gamma_k\right|^2}= |\lambda|\sum_{k=1}^{\infty}\frac{\frac{\mathcal{A}}{k^{\alpha}}+\varphi(k)}{\left|\lambda\right|^2+\left|\mathcal{B}k^{\beta}+\phi(k)\right|^2}.
\end{align} As result of estimating of the modulus of the difference
 \begin{align*}
\left|\frac{\left[\frac{\mathcal{A}}{k^{\alpha}}+\varphi(k)\right] |\lambda|}{|\lambda|^2+(\mathcal{B}k^{\beta}+\phi (k))^2}-\frac{\frac{\mathcal{A}}{k^{\alpha}}|\lambda|}{|\lambda|^2+\mathcal{B}^2k^{2\beta}}\right|&\leq \frac{\varphi(k)\left(|\lambda|^3+|\lambda|\mathcal{B}^2k^{2\beta}\right)+\frac{\mathcal{A}|\lambda|}{k^{\alpha}}\left(2\mathcal{B}k^{\beta}\phi(k)+\phi^2(k)\right)}{(|\lambda|^2+(\mathcal{B}k^{\beta}+\phi (k))^2)(|\lambda|^2+\mathcal{B}^2k^{2\beta})}\\
&=\frac{O\left(\frac{1}{k^{4\beta}}\right)+\frac{1}{|\lambda|^2}\left[O\left(\frac{1}{k^{2\beta}}\right)+O\left(\frac{1}{k^{2\beta}}\right)+O\left(\frac{1}{k^{2\beta+1}}\right)\right]}{|\lambda|k^{\alpha+1}\left(\frac{1}{k^{2\beta}}+\frac{\left(\mathcal{B}+O\left(\frac{1}{k}\right)\right)^2}{|\lambda|^2}\right)\left(\frac{1}{k^{2\beta}}+\frac{\mathcal{B}^2}{|\lambda|^2}\right)}
\end{align*} the following estimate  is established immediately
\begin{align}
\left|\frac{\left[\frac{\mathcal{A}}{k^{\alpha}}+\varphi(k)\right] |\lambda|}{|\lambda|^2+(\mathcal{B}k^{\beta}+\phi (k))^2}-\frac{\frac{\mathcal{A}}{k^{\alpha}}|\lambda|}{|\lambda|^2+\mathcal{B}^2k^{2\beta}}\right| \lesssim \frac{1}{|\lambda|k^{\alpha+1}}. \label{approximation-of-derivate-of-kernel-and-its-part-without-varphi-function}
\end{align}
Hence, from the estimates \eqref{estimate-of-derivate-kernel-function-by-lambda} and \eqref{approximation-of-derivate-of-kernel-and-its-part-without-varphi-function} and making the change of variables $\mathcal{B}x^{\beta}=|\lambda|t$ the chain of inequalities is attained.
 \begin{align*}
\left|\lambda \widehat{K}^{'}(\lambda)\right|&\lesssim \left|\sum_{k=1}^{\infty}\frac{k^{-\alpha}\mathcal{A}|\lambda|}{(|\lambda|^2+\mathcal{B}^2k^{2\beta})}\right|+\frac{1}{|\lambda|}\sum_{k=1}^{\infty}\frac{1}{k^{\alpha+1}}
\leq \frac{\mathcal{A}|\lambda|}{|\lambda|^2+\mathcal{B}^2}+\left|\int_{2}^{\infty}\frac{x^{-\alpha}\mathcal{A}|\lambda|dx}{(|\lambda|^2+\mathcal{B}^2x^{2\beta})}\right|+\frac{1}{|\lambda|}\sum_{k=1}^{\infty}\frac{1}{k^{\alpha+1}}\\
&\lesssim \frac{\mathcal{A}}{|\lambda|}+\frac{\mathcal{A}\mathcal{B}^{r-1}}{|\lambda|^r}\left|\int_{\frac{\mathcal{B}2^{\beta}}{|\lambda|}}^{\infty}\frac{dt}{t^{r}(1+t^2)}\right|+\frac{1}{|\lambda|}\sum_{k=1}^{\infty}\frac{1}{k^{\alpha+1}}.
\end{align*} Thus, $\left|\lambda \widehat{K}^{'}(\lambda)\right| \to 0$ when $|\lambda| \to +\infty$.  Lemma \ref{lemma-argument-pi-delta} is proved.

\section*{Acknowledgements}

This ~work was ~supported ~by the ~Mexican ~Center ~for ~Economic ~and ~Social ~Studies (CEMEES) and the Russian Science Foundation project  no. 17-11-01215. 

The authors are very grateful to Professors A. A. Shkalikov and A. L. Skubachevskii for their advices and constructive comments.


\begin{thebibliography}{99} 
\bibitem{GMJ} G. Amendola, M. Fabrizio, J. M. Golden, \emph{Thermodynamics of Materials
with Memory}, Theory and applications. Springer, New York, Dordrecht, Heidelberg, London, 2012.
\bibitem{Dafermos} C. M. Dafermos, Asymptotic stability in viscoelasticity, \emph{Archive for Rational Mechanics and Analysis}, \textbf{37}, 1970, 297--308.
\bibitem{MFBL} M. Fabrizio, B. Lazzari, On the existence and the asymptotic stability of solutions for linearly viscoelastic solids,  \emph{Archive for Rational Mechanics and Analysis}, \textbf{116}, 1991, 139--152.
\bibitem{PG} M. E. Gurtin, A. C. Pipkin, A general theory of heat conduction with finite wave speeds,  \emph{Archive for Rational Mechanics and Analysis},  \textbf{31}, 1968, 113--126.
\bibitem{GK} R. A. Guyer, J. A. Krumhansl, Solution of the linearized phonon Boltzmann equation, \emph{Phys. Rev.}  \textbf{148}, 1966, 766--778.
\bibitem{PI} S. Ivanov, L. Pandolfi, Heat equations with memory: lack of controllability to rest,  \emph{Journal of Mathematical Analysis and Applications,} \textbf{355}, 2009, 1--11.
\bibitem{LM} J. L. Lions, E. Magenes, \emph{Probleme aux limites non homog\`enes et applications} (in French)  [Non-Homogeneous Boundary  Problems and its Applications],  \textbf{1}, S. A. Dunod, Paris, 1968; Mir, Moscow, 1971.
\bibitem{MKR} R. K. Miller, Volterra integral equations in a Banach space, \emph{Funkcial. Ekvac.}, \textbf{8}, 1975, 163--193.
\bibitem{MKR1} R. K. Miller, An integrodifferential equation for rigid heat conductors with memory,  \emph{J. Math. Anal. Appl.}, \textbf{66}, 1978, 313--332.
\bibitem{MRWR}R. K. Miller, R. L. Wheeler, Well-posedness and stability of linear Volterra integrodifferential equations in abstract spaces, \emph{Funkcial. Ekvac.}, \textbf{21}, 1978, 279--305.
\bibitem{MRWD} R. K. Miller, W. Desch, Exponential stabilization of Volterra integrodifferential equations in Hilbert space, \emph{Journal of Differential Equations}, \textbf{70}, 1987, 366--389.
\bibitem{JMF} J. E. Mu\~noz Rivera, M. Grazia Naso, F. M. Vegni, Asymptotic behavior of the energy for a class of a weakly dissipative second-order systems with memory,  \emph{Journal of Mathematical Analysis and Applications}, \textbf{286}, 2003, 692--704.
\bibitem{JME} J. E. Mu\~noz Rivera, M. Grazia Naso, E. Vuk, Asymptotic behavior of the energy for electromagnetic systems with memory, \emph{Mathematical Methods in the Applied Sciences}, \textbf{27}, 2004, 819--841.
\bibitem{JM} J. E. Mu\~noz Rivera, M. Grazia Naso, On the decay of the energy for systems with memory and indefinite dissipation, \emph{Asymptotic Analysis}, \textbf{49}, 2006, 189--204.
\bibitem{P} L. Pandolfi, The controllability of the Gurtin--\allowbreak Pipkin equations: a cosine operator approach, \emph{Applied Mathematics and Optimization}, \textbf{52}, 2005, 143--165.
\bibitem{RomeoRautian} R. Perez Ortiz, N. A. Rautian,  Representation of solutions of integro-differential equations with kernels depending on the parameter, \emph{Differential Equations}, \textbf{53} (1), 2017, 139--143.
\bibitem{RPV2} R. Perez Ortiz, V. V. Vlasov,  Spectral analysis of integrodifferential equations arising in the theory of viscoelasticity and thermal physics, \emph{Mathematical Notes}, \textbf{98} (4), 2015, 689--693.
\bibitem{RP1} R. Perez Ortiz, Spektralnii analiz integrodiferentsialnyx urabnenyi c yadrami zavicyashimi ot parametra (in russian) [Spectral Analysis of integrodifferential equations with kernels depending on parameter], in \emph{Proceedings of Moscow Institute of Physics and Technology (MIPT)}, \textbf{7} (2), 2015, 27--38.
\bibitem{VR} V. V. Vlasov, N. A. Rautian, Well-defined solvability and spectral analysis of abstract hyperbolic integrodifferential equations,  \emph{Journal of Mathematical Sciences}, \textbf{179}, 2011, 390--414.
\bibitem{VRShamaev}  V. V. Vlasov, N. A. Rautian, A. S. Shamaev, Spectral analysis and correct solvability of abstract integrodifferential equations arising in thermophysics and acoustics, \emph{Journal of Mathematical Sciences},  \textbf{190}, 2013, 34--65.
\bibitem{VR1}  V. V. Vlasov, N. A. Rautian, Spectral analysis of hyperbolic Volterra integrodifferential equations, \emph{Doklady mathematics}, \textbf{92} (2), 2015, 590--593.
\bibitem{VR2} V. V. Vlasov, N. A. Rautian, Spectral analysis and representations of solutions of abstract integrodifferential equations in Hilbert space, \emph{Operator Theory: Advances and Applications}, \textbf{236}, 2014, 517--535.
\bibitem{ESP} E. Sanchez-Palencia, Nonhomogeneous media and vibration theory, Lecture Notes in physics, Springer-Verlag, Berlin Heidelberg, New York, 1980.
\bibitem{Vegni} F. M. Vegni, Dissipativity of a conserved phase-field system with memory, \emph{Discrete and continuous dynamical systems}, \textbf{9}, 2003, 949--968.
 \bibitem{shapiro} J. H. Shapiro,  Composition operators and classical function theory. New York: Springer, 1993.
 \bibitem{shapiroBourdon} J. H. Shapiro, P. S. Bourdon, Cyclic phenomena for composition operators. \emph{Memoirs of the American Mathematical Society}, Vol. 125, No. 596, 1997.
  \bibitem{PORVVV} V. V. Vlasov, Perez Ortiz R., Correct solvability of Volterra integrodifferential equations in Hilbert space. \emph{Electronic Journal of Qualitative Theory of Differential Equations (EJQTDE)}, 2016, No. 31, 1--17.
 \bibitem{Vollunaranovich} L. I. Volkoviskii, G. L. Lunts, I. G. Aramovich, Sbornik sadach po teorii funksii kompleksnovo peremenovo (in russian) [Compilation of problems on the theory of functions of a complex variable]. Uchebnoe posobye, izdanie 4, Fizmatlit, 2002, 312 pp.
\bibitem{vlasovmrau2016} V. V. Vlasov, N. A. Rautian, Spektralnii analisis funktsionalno-differentsialnyx urabnenii (in russian)[Spectral analysis of functional differential equations], Monograph--M. editorial MAKS Press, Moscow,  2016, 488 pp.
\end{thebibliography}
\end{document}